\documentclass{siamart171218}
\usepackage{color}
\usepackage{latexsym}
\usepackage{dsfont}
\usepackage{amssymb}
\usepackage{amsmath, amsfonts,amssymb,theorem,euscript,array,enumerate,amsfonts,mathrsfs}

\usepackage{hyperref}
\usepackage{enumitem}
\setlist{leftmargin=*, topsep=0.5em, parsep=0pt, itemsep=1em, labelindent=0pt, align=left}
\usepackage{nicefrac}

\headers{Neural networks algorithms for stochastic control}{C. Hur\'e, H. Pham, A. Bachouch and N. Langren\'e}

\newsiamremark{remark}{Remark}

\usepackage{accents}
\newlength{\dhatheight}

\newcommand{\ep}{	\begin{flushright}
		$\blacksquare$
\end{flushright}}
\def \Prod{\displaystyle\prod}

\def \Inf{\displaystyle\inf}
\def \Sup{\displaystyle\sup}

\def \Max{\displaystyle\max}
\def \Min{\displaystyle\min}
\def \b1{\bf{1}}

\def \A{\mathbb{A}}

\def \N{\mathbb{N}}
\def \R{\mathbb{R}}

\def \E{\mathbb{E}}
\def \F{\mathbb{F}}

\def \P{\mathbb{P}}

\def\esssup_#1{\underset{#1}{\mathrm{ess\,sup\, }}}

\def\argmin_#1{\underset{#1}{\mathrm{argmin\, }}}

\def \Ac{{\cal A}}
\def \Bc{{\cal B}}
\def \Cc{{\cal C}}

\def \Fc{{\cal F}}

\def \Oc{{\cal O}}

\def \Xc{{\cal X}}
\def \Vc{{\cal V}}

\def \Op{{\cal O}_{\P}}

\def \eps{\varepsilon}
\newcommand*\diff{\mathop{}\!\mathrm{d}}

\def \ep{\hbox{ }\hfill$\Box$}

\def\reff#1{{\rm(\ref{#1})}}

\def\beqs{\begin{eqnarray*}}
\def\enqs{\end{eqnarray*}}
\def\beq{\begin{eqnarray}}
\def\enq{\end{eqnarray}}

\makeatletter
\newcommand{\setword}[2]{%
	\phantomsection
	#1\def\@currentlabel{\unexpanded{#1}}\label{#2}%
}
\makeatother

\ifpdf
\hypersetup{
	pdftitle={Deep neural networks algorithms for stochastic control problems on finite horizon: convergence analysis},
	pdfauthor={C\^ome Hur\'e, Huy\^en Pham, Achref Bachouch and Nicolas Langren\'e}
}
\fi

\begin{document}
 
\title{
Deep neural networks algorithms for\\ stochastic control problems 
on finite horizon: convergence analysis 
}

 \author{C\^ome \textsc{Hur\'e}\thanks{LPSM, Universit\'e de Paris (Paris Diderot)
 		(\email{hure at lpsm.paris}).}
 	\and Huy\^en \textsc{Pham}\thanks{LPSM, Universit\'e de   Paris (Paris Diderot) and CREST-ENSAE
 		(\email{pham at lpsm.paris}). The work of this author is su\-pported by the ANR project CAESARS (ANR-15-CE05-0024), and also by FiME and the ``Finance and Sustainable Development'' EDF - CACIB Chair.}
 	\and Achref \textsc{Bachouch}\thanks{Department of Mathematics, University of Oslo, Norway.
	(\email{achrefb at math.uio.no}). The  author's research was carried out with support of the Norwegian Research Council, within the research project Challenges in Stochastic Control, Information and Applications (STOCONINF), project number 250768/F20.}
 	\and Nicolas \textsc{Langren\'e}\thanks{CSIRO Data61, Melbourne, Australia
	(\email{nicolas.langrene at csiro.au}).}
 	}

\maketitle 

\begin{abstract}
This paper develops  algorithms for high-dimensional stochastic control pro\-blems based on deep learning and dynamic programming. Unlike classical approximate dynamic programming approaches, we first approximate the optimal policy  by means of neural networks in the spirit of deep reinforcement learning, and then the value function by Monte Carlo regression. This is achieved in the dynamic programming recursion by performance or hybrid iteration, and regress now methods from numerical probabilities.  We provide a theoretical justification of these algorithms. 
Consistency and rate of convergence for the control and value function estimates are analyzed and expressed in terms of the universal approximation 
error of the neural networks, and of the statistical error when estimating network function, leaving aside the optimization error.  
Numerical results on various applications are presented in a companion paper \cite{bacetal18} and illustrate the performance of the proposed algorithms. 
\end{abstract} 

%\vspace{4mm}


%\vspace{4mm}


\section{Introduction} \label{secintro}
A large class of dynamic decision-making problems under uncertainty can be mathematically modeled as discrete-time stochastic optimal control problems in finite horizon. This paper is devoted to the analysis of novel probabilistic numerical algorithms based on neural networks for solving such problems. Let us consider the following discrete-time stochastic  control problem over a finite horizon $N$ $\in$ $\N\setminus\{0\}$. 
 %over a finite horizon $N$ $\in$ $\N\setminus\{0\}$. 
 The dynamics of the controlled state process $X^\alpha$ $=$ $(X^\alpha_n)_{n}$ valued in $\Xc$ $\subset$ $\R^d$  is given by
 \beq \label{dynX}
 X_{n+1}^\alpha &=& F(X_n^\alpha,\alpha_n,\eps_{n+1}), \;\;\; n=0,\ldots,N-1, \; X_0^\alpha = x_0 \in \R^d, 
 \enq
 where $(\eps_n)_{n}$  is a sequence of i.i.d. random variables valued in some Borel space $(E,\Bc(E))$, and defined on some probability space $(\Omega,\Fc,\P)$ equipped with the filtration   $\F$ $=$ $(\Fc_n)_n$ generated by the noise $(\eps_n)_n$ ($\Fc_0$ is the trivial $\sigma$-algebra), the control $\alpha$ $=$ $(\alpha_n)_{n}$ is an 
 $\F$-adapted process valued in  $\A$ $\subset$ $\R^q$, and $F$ is a measurable function from $\R^d\times\R^q\times E$ into $\R^d$. 
   
Given a running cost function $f$ defined on $\R^d\times\R^q$,  a terminal cost function $g$ defined on $\R^d$, 
%and a discount factor $\rho$ $\in$ $[0,1]$, 
the cost functional associated with a control process $\alpha$ is
$
 J(\alpha) = \E \left[ \sum_{n=0}^{N-1}  f(X_n^\alpha,\alpha_n) +  g(X_N^\alpha) \right]. 
 $
 The set  $\Cc$ of admissible control is the set of control processes $\alpha$ satisfying some integrability conditions  ensuring  that the cost functional $J(\alpha)$ is well-defined and finite. The control problem, also called Markov decision process (MDP), is formulated as  
\beq \label{defcontrol}
V_0(x_0) & := & \inf_{\alpha\in\Cc} J(\alpha),
\enq
and the goal is to find  an optimal control $\alpha^*$ $\in$ $\Cc$, i.e.,  attaining the optimal value: $V_0(x_0)$ $=$ $J(\alpha^*)$. 
Notice that problem  \reff{dynX}-\reff{defcontrol} may also be viewed as the time discretization of a continuous time stochastic control problem, in which case $F$ is typically the Euler scheme for a controlled diffusion process, and $V_0$ is the discrete-time approximation of a fully nonlinear Hamilton-Jacobi-Bellman equation.

Problem \reff{defcontrol} is tackled by the dynamic programming approach, and we introduce the standard notations for MDP:  denote by 
$\{P^a(x,dx')$, $a$ $\in$ $\A$, $x$ $\in$ $\Xc \}$, the family of  transition probabilities associated with the controlled (homogenous) Markov chain \reff{dynX},  given by $ P^a(x,dx') :=  \P\left[ F(x,a,\eps_{1}) \in dx' \right] $
and for any measurable function $\varphi$ on $\Xc$:
$
P^a \varphi(x) = \int \varphi(x') P^a(x,dx')   =  \E \left[ \varphi\big(F(x,a,\eps_{1}) \big) \right] . 
$
With these notations, we have  for any measurable function $\varphi$ on $\Xc$, for any $\alpha$ $\in$ $\Cc$, 
$
\E [ \varphi(X_{n+1}^\alpha) | \Fc_n ] = P^{\alpha_n} \varphi(X_n^\alpha), \;\;\; \forall \; n \in \N.  
$
The optimal value $V_0(x_0)$ is then determined in backward induction starting from the terminal condition 
$
V_N(x) =g(x), \; x \in \Xc,  
$
and  by the dynamic programming (DP) formula, for $n$ $=$ $N-1,\ldots,0$: 
\begin{equation} \label{DP}
\left\{
\begin{array}{rcl}
Q_n(x,a) &=& f(x,a) +   P^a V_{n+1}(x), \;\;\; x \in \Xc, \; a \in \A, \\
V_n(x) &=& \Inf_{a\in \A} Q_n(x,a),
\end{array}
\right.
\end{equation}
The function $Q_n$ is called optimal state-action value function, and $V_n$ is the (optimal) value function. 
Moreover, when the infimum is attained in the DP formula at any time $n$ by $a_n^*(x)$, we get an optimal control in feedback form given by: $\alpha^*$ $=$ 
$(a_n^*(X_n^*))_n$ where $X^*$ $=$ $X^{\alpha^*}$  is the Markov process defined by
\beqs
X_{n+1}^* &=& F(X_n^*,a_n^*(X_n^*),\eps_{n+1}), \;\;\;  n=0,\ldots,N-1, \;\; X_0^* = x_0. 
\enqs

The implementation of the DP formula requires the knowledge and explicit computation of the transition probabilities $P^a(x,dx')$. 
%as well as  of the reward function $f_n(x,a)$. 
In situations when they are unknown, this leads to the problematic of reinforcement lear\-ning for computing the optimal control and value function by relying on simulations of  the environment.  The challenging tasks from a numerical point of view  are then twofold:
\begin{itemize}
\item[1.] {\it Transition probability operator.} Calculations  for any  $x$ $\in$ $\Xc$, and any $a$ $\in$ $\A$, of  $P^a V_{n+1}(x)$,  for $n$ $=$ $0,\ldots,N-1$.     
This is  a computational challenge in high dimension $d$  for the state space with the ``curse of dimensionality"  due to the explosion of grid points in deterministic methods. 
\item[ 2.]  {\it Optimal control.} Computation of  the infimum in $a$ $\in$ $\A$ of $ f(x,a) + P^a V_{n+1}(x)$ for fixed $x$ and $n$, 
and of $\hat a_n(x)$ attaining the minimum if it exists.  This is also a computational challenge especially in high dimension $q$ for the control space. 
\end{itemize}

The classical probabilistic numerical methods based on DP for solving the MDP are sometimes 
called approximate dynamic programming methods, see e.g. \cite{bertso96}, \cite{pow11}, and consist basically of the two following steps: 
\begin{itemize}
\item[(i)] Approximate the $Q_n$. This can be performed by regression Monte Carlo (RMC) techniques or quantization.  RMC is  typically done by least squares linear regression on a set of basis function following the popular approach by Longstaff and Schwartz \cite{LS01} initiated for the Bermudan option problem, where the suitable choice of basis functions might be delicate. Conditional expectation can also be approximated by regression on neural networks as in \cite{kohkrztod10} for the American option problem and  appears as a promising and efficient alternative in high dimension.  
\item[(ii)]  Control search. Once we get an approximation $(x,a)$ $\mapsto$ $\hat Q_n(x,a)$ of the $Q_n$ value function, the optimal control $\hat a_n(x)$ which achieves the minimum over $a$ $\in$ $\A$ of $Q_n(x,a)$ can be obtained either by an exhaustive search when $\A$ is discrete (with relatively small cardinality), or by a (deterministic) gradient-based algorithm for continuous control space (with relatively small dimension). 
\end{itemize}

Recently,  numerical methods by direct approximation, without DP, have been deve\-loped and made implementable thanks to the power of computers: the basic idea is to focus directly on the control approximation by considering feedback control (policy) in a parametric form:
$
a_n(x) = A(x;\theta_n), \;\;\; n=0,\ldots,N-1,
$
for some given function $A(.,\theta_n)$ with parameters $\theta$ $=$ $(\theta_0,\ldots,\theta_{N-1})$ $\in$  $\R^{q\times N}$, and minimize over $\theta$ the parametric functional
$
\tilde J(\theta) =\E \left[ \sum_{n=0}^{N-1}  f(X_n^A,A(x;\theta_n)) +  g(X_N^A) \right], 
$
where $(X_n^A)_n$ denotes the controlled process with feedback control $(A(.,\theta_n))_n$. This approach was first adopted in \cite{kouetal16}, which used the EM algorithm for optimizing over the parameter $\theta$,  and was further investigated in \cite{hanE16}, \cite{Ehanjen17}, \cite{PHL17}, which considered deep neural networks (DNN) for the parametric feedback control, and stochastic gradient descent methods (SGD) for computing the optimal parameter 
$\theta$. The theoretical foundation of these DNN algorithms has been recently investigated in \cite{hanlon18}. 
Let us mention that DNN approximation in stochastic control has already been explored in the context of reinforcement learning (RL) (see \cite{bertso96} and \cite{sutbar98}), and is called deep reinforcement learning in the artificial intelligence community \cite{mni15} (see also  \cite{li17} for a recent survey) but usually for infinite horizon (stationary) control problems.

In this paper, we combine different ideas from the mathematics (numerical probability) and the computer science (reinforcement learning) communities to propose and compare several algorithms based on dynamic programming (DP), and deep neural networks (DNN) for the approximation/learning of (i) the optimal policy, and then of (ii) the value function.  Notice that this differs from the classical approach in DP recalled above, where we first approximate the $Q$-optimal state/control value function, and then approximate the optimal control. Our learning of the 
optimal policy is achieved in the spirit of \cite{hanE16} by DNN, but sequentially in time through DP instead of a global learning over the whole period $0,\ldots,N-1$. Once we get an approximation of the optimal policy, we approximate the value function by Monte Carlo (MC) regression based on simulations of the forward process with the approximated optimal control. In particular, we avoid  the issue of {\it a priori} endogenous simulation of the controlled process in the classical $Q$-approach. The MC regressions for the approximation of the optimal policy and/or value function, rely  on performance iteration (PI) or hybrid iteration (HI).  Numerical results on several applications are devoted to a companion paper \cite{bacetal18}. 
The theoretical contribution of the current  paper is to provide a detailed convergence analysis of our  two proposed algorithms:  Theorem \ref{theo::algoIII}  for the {\it NNContPI } algorithm based on control learning by performance iteration with DNN,  and
Theorem \ref{theoKohler} for the {\it Hybrid-Now} algorithm based on control learning by DNN and then value function learning by the regress-now method. We rely mainly on arguments from statistical learning and non parametric regression as developed notably in the book \cite{gyokoh02} to give estimates of  approximated control and value function in terms of the universal approximation error of the  neural networks, and of  the statistical error in the estimation of network functions. In particular, the estimation of the optimal control using the control that minimizes the empirical loss is studied based on the uniform law of large numbers (see \cite{gyokoh02}); and the estimation of the value function in Hybrid-Now relies on Lemma \ref{lem:H1}, proved in \cite{kohkrztod10}, which bounds the difference between the conditional expectation and the best approximation of the empirical conditional expectation by neural networks.
In practical implementation, the computation of the approximate optimal policy in the minimization of the empirical loss function is performed by a stochastic gradient descent algorithm, 
which leads to the so-called   {\it optimization error}, see e.g. \cite{berfor13}, \cite{jenetal18}.  
 We do not address in the current paper  this  type of error in machine learning algorithms, and we refer to the recent paper \cite{becjenkuc20} for a full error analysis of deep neural networks training. 

The paper is organized as follows.  We recall in Section \ref{secpreli} some basic results about deep neural networks (DNN)  and stochastic optimization gradient descent methods. Section \ref{sec::listeAlgos} is devoted to the description of our two algorithms. We  analyze in Section \ref{secconv} the convergence of the algorithms. Finally the Appendices collect some Lemmas used in the proof of the convergence results. 

\section{Preliminaries on DNN and SGD} \label{secpreli}

\subsection{Neural network approximations} \label{secNN}

Deep Neural networks (DNN) aim to approximate (complex nonlinear)  
functions defined on finite-dimensional space, and in contrast with the usual additive approximation theory built via basis functions, 
like polynomial, they rely on composition of layers of simple functions.  The relevance of neural networks  comes from the universal approximation theorem and the 
Kolmogorov-Arnold representation theorem (see \cite{kol61}, \cite{cyb89} or \cite{hor91}), and  this has been shown to be successful in numerous practical applications. 

We consider here a feedforward artificial network (also called multilayer perceptron) for the approximation of the optimal policy (valued in $\A$ $\subset$ $\R^q$) 
and the value function (valued in $\R$), both defined on the state space $\Xc$ $\subset$ $\R^d$. The architecture is represented by  functions    
$
x  \; \in \;  \Xc  \longmapsto  \Phi(z;\theta) \; \in \R^o, 
$
with $o$ $=$ $q$ or $1$ in our context, and where $\theta$ $\in$ $\Theta$ $\subset$ $\R^p$ are the weights (or parameters) of the neural networks.

\vspace{3mm}

\subsection{Stochastic optimization in DNN} \label{secoptimsto}

Approximation by means of DNN  requires a stochastic optimization with respect to a set of parameters, which can be written in a generic form as
\beq \label{Loss}
\inf_\theta \E \big[ L(Z;\theta) \big],
\enq
where $Z$ is a random variable from which the training samples $Z^{(m)}$, $m$ $=$ $1,\ldots,M$ are drawn, and $L$ is a loss  function involving DNN with parameters $\theta$ $\in$ 
$\R^p$, and typically differentiable w.r.t. $\theta$ with known gradient $D_\theta L(Z^{(m)})$.  

Several algorithms based on stochastic gradient-descent are often used in practice for the search of infimum in \reff{Loss}; see e.g.  basic SGD, AdaGrad, RMSProp, or Adam.

In the theoretical study of the algorithms proposed in this paper, we did not address the convergence of the gradient-descent algorithms,  
%(as done e.g. in \cite{jenetal18}), 
and rather assumed that the agent successfully solved the minimization problem \eqref{Loss}, i.e. minimized the empirical loss $\sum_{m=1}^{M}L(Z^{(m)})$ using (non stochastic) gradient-descent, a stochastic gradient-descent method with many epochs or any other methods that allow overfitting. One can then rely on statistical learning results as in \cite{gyokoh02} and \cite{kohkrztod10} that are specific to neural networks. Notice that overfitting should be avoided in practice, and so our framework remains theoretical.

\section{Description of the algorithms}
\label{sec::listeAlgos}

Let us introduce a set $\Ac$ of neural networks for approximating optimal policies, that is a set of parametric functions $x$ $\in$ $\Xc$ $\mapsto$ $A(x;\beta)$ $\in$ $\A$, with parameters 
$\beta$ $\in$ $\R^l$, and a set $\Vc$ of neural networks functions for approximating value functions, that is a set of parametric functions 
$x$ $\in$ $\Xc$ $\mapsto$ $\Phi(x;\theta)$ $\in$ $\R$, with parameters $\theta$ $\in$ $\R^p$. 

We are also given at each time $n$ a probability measure $\mu_n$ on the state space $\Xc$, which we refer to as a training distribution. 
Some comments about the choice of the training measure are discussed in Section \ref{sectraining}.

\subsection{Control learning by performance iteration}  \label{secAlgoPI}
 
This algorithm, referred in short as {\it NNcontPI} algorithm, is designed as follows: 

For $n$ $=$ $N-1,\ldots,0$, keep track of the approximated optimal policies $\hat a_k$, $k$ $=$ $n+1,\ldots,N-1$, and
 approximate the optimal policy at time $n$ by 
$\hat a_n$  $=$ $A(.;\hat\beta_n)$ with 
\beq 
 \hat\beta_n  & \in & \argmin_{\beta \in \R^l}  \E \left[ f(X_n,A(X_n;\beta)) + \sum_{k=n+1}^{N-1} f(\hat X_k^\beta,\hat a_k(\hat X_k^\beta)) + g(\hat X_N^\beta) \right],  \label{algoPI}
\enq
where $X_n$ $\sim$ $\mu_n$,  $\hat X_{n+1}^\beta$ $=$ $F(X_n,A(X_n;\beta),\eps_{n+1})$ $\sim$ $P^{A(X_n;\beta)}(X_n,dx')$, and for $k$ $=$ $n+1,\ldots,N-1$, 
$\hat X_{k+1}^\beta$ $=$ $F(\hat X_{k}^\beta,\hat a_k(\hat X_k^\beta),\eps_{k+1})$ $\sim$ $P^{\hat a_k(\hat X_k^\beta)}(\hat X_k^\beta,dx')$. 
Given estimate $\hat a_k^M$ of  $\hat a_k$, $k$ $=$ $n+1,\ldots,N-1$,  the approximated policy $\hat a_n$ is estimated by using a training sample 
$\left(X_n^{(m)},(\eps_{k+1}^{(m)})_{k=n}^{k=N-1}\right)$, $m$ $=$ $1,\ldots,M$ of  $\left(X_n,(\eps_{k+1})_{k=n}^{k=N-1}\right)$ for simulating $\left(X_n,(\hat X_{k+1}^\beta)_{k=n}^{k=N-1}\right)$, 
and optimizing $\beta$ of $A(.;\beta)$ $\in$ $\Ac$, by SGD (or its variants) as described in Section \reff{secoptimsto}. 

%\noindent $\tri$ We then get an estimate of the optimal policy at any time $n$ $=$ $0,\ldots,N-1$ by:
%\beqs
%\hat a_n^M & = & A(.;\hat\beta_n^M)  \; \in \; \Ac, 
%\enqs
%where $\hat\beta_n^M$ is the estimated optimal parameter resulting from the SGD in  \reff{algoPI} with a training sample of size $M$.   
The estimated value function at time $n$  is then given as
\begin{equation} \label{estimValuePI}
\hat V_n^M(x) =
%\frac{1}{M} \sum_{m=1}^M  \Big[  \sum_{k=n}^{N-1}  f(\hat X_k^{(m)},\hat a_k^M(\hat X_k^{(m)})) +  g(\hat X_N^{(m)}) \Big], 
\E_M \left[ \sum_{k=n}^{N-1}  f(\hat X_k^{n,x},\hat a_k^M(\hat X_k^{n,x})) +  g(\hat X_N^{n,x}) \right], 
\end{equation}
where $\E_M$ is the expectation conditioned on the training set used for computing $\left(\hat a_k^M\right)_k$; and 
$\left(\hat X_k^{n,x}\right)_{k=n}^{N}$ is driven by the estimated optimal controls. The dependence of the estimated value function $\hat V_n^M$ upon the training samples $X_k^{(m)}$, for $m$ $=$ $1,\ldots,M$, used at time $k$ $=$ $n,\ldots,N$, 
is emphasized through the exponent $M$ in the notations.   
 
\vspace{2mm}

\begin{remark}
{\rm The {\it NNcontPI} algorithm can be viewed as a combination of the DNN algorithm designed in \cite{hanE16} and dynamic programming. 
In the algorithm presented in \cite{hanE16}, which totally ignores the dynamic programming principle, 
one learns all the optimal controls $A(.;\beta_n)$, $n$ $=$ $0,\ldots,N-1$ at the same time, by performing one unique stochastic gradient descent. 
This is efficient as all the parameters of all the NN are getting trained at the same time, using the same mini-batches. However, when the number of layers of the global neural network 
gathering all the $A(.;\beta_n)$, for $n=0,\ldots,N-1$ is large (say $ \sum_{n=0}^{N-1} \ell_n$ $\geq$ $100$, where $\ell_n$ is the number of layers in $A(. ,\beta_n)$), then one is likely 
to observe vanishing or exploding gradient problems that will affect the training of the weights and bias of the first layers of the global NN (see  \cite{ger17}  for more details). 
%Therefore, it may be more reasonable to make use of the dynamic programming structure when $N$ is large, and learn the optimal policy sequentially as proposed in our {\it NNcontPI} Algo. 
%Notice that a similar idea was already used in \cite{guylab10} in the context of uncertain volatility model where the authors  use a specific parametrization 
%for the feedback control instead of a DNN adopted more generally here.
}
\ep
\end{remark}

\begin{remark}
{\rm The {\it NNcontPI} algorithm does not require value function iteration, but instead is based on 
performance iteration by keeping track of the estimated optimal policies computed in backward recursion.  The value function is then computed in \reff{estimValuePI} as the gain functional 
associated with the estimated optimal policies $(\hat a_k^M)_k$. Consequently, it provides usually a low bias estimate but induces 
possibly high variance estimate and large complexity, especially when $N$ is large.
 }
\ep
\end{remark}

\subsection{Control learning by hybrid iteration}

Instead of keeping track of all the approximated optimal policies as in the {\it NNcontPI} algorithm, 
we use an approximation of the value function at time $n+1$ in order to compute the optimal policy at time $n$. The approximated value function is then updated at time $n$. This leads to the following algorithm: 

\noindent \rule{\linewidth}{.5pt}

\noindent  {\bf Hybrid-Now algorithm}  

\begin{itemize}
\item[1.] {\it Initialization}: $\hat V_N$ $=$ $g$ 
\item[2.]   For $n$ $=$ $N-1,\ldots,0$, 
\begin{itemize}
\item[(i)] Approximate the optimal policy at time $n$ by $\hat a_n$ $=$ $A(.;\hat\beta_n)$ with 
\beq \label{algoHI}
\hat\beta_n & \in & \argmin_{\beta \in \R^l}  \E \left[ f(X_n,A(X_n;\beta)) + \hat V_{n+1}(X_{n+1}^{A(.,\beta)}) \right], 
\enq
where $X_n$ $\sim$ $\mu_n$,  $\hat X_{n+1}^{A(.,\beta)}$ $=$ $F(X_n,A(X_n;\beta),\eps_{n+1})$ $\sim$ $P^{A(X_n;\beta)}(X_n,dx')$. 
\item[(ii)] {\it Updating}: approximate the value function by NN:
\beq \label{algoHItheta}
\hat V_n^M & \in & \argmin_{\Phi(.;\theta) \in \Vc} \E \big| f(X_n,\hat a_n^M(X_n)) + \hat V_{n+1}^M(X_{n+1}^{\hat a_n^M}) \; - \; \Phi(X_n;\theta) \big|^2
\enq
using samples $X_n^{(m)}$, $X_{n+1}^{\hat a_n^M,(m)}$, $m$ $=$ $1,\ldots,M$ of $X_n$ $\sim$ $\mu_n$, and $X_{n+1}^{\hat a_n^M,(m)}$  of $X_{n+1}^{\hat a_n^M}$
\end{itemize}
\end{itemize}

\noindent \rule{\linewidth}{.5pt}

\vspace{2mm}

The approximated policy $\hat a_n$ is estimated by using a training sample $(X_n^{(m)},\eps_{n+1}^{(m)})$, $m$ $=$ $1,\ldots,M$ of $(X_n^{},\eps_{n+1}^{})$ to simulate 
$\left(X_n,X_{n+1}^{A(.;\beta)}\right)$, 
and optimizing $A(.;\beta)$ $\in$ $\Ac$ over the parameters $\beta$ $\in$ $\R^l$ and the expectation in \reff{algoHI} by stochastic gradient descent methods as described in Section \reff{secoptimsto}. We then get an estimate $\hat a_n^M$ $=$ $A\left(.;\hat\beta_n^M\right)$. 
The  value function is represented by neural network and its optimal parameter is estimated by MC-regression in \reff{algoHItheta}.

\subsection{Training sets  design}
\label{sectraining}
We discuss here the choice of the training measure $\mu_n$ used to generate the training sets on which will be computed the estimates.
Two cases are considered in this section. The first one is a knowledge-based selection, relevant when the controller knows with a certain degree of confidence where the process has to be driven in order to optimize her cost functional.
The second case, on the other hand, is when the controller has no idea where or how to drive the process to optimize the cost functional. 
Separating the construction of the training set from the training procedure of the optimal control and/or value function simplifies the theoretical study of our algorithms. It is nevertheless common to merge them in practice, as discussed in Remark \ref{rmk:merge_trainings}. 

\vspace{3mm}

\noindent {\bf Exploitation only strategy} 
In the knowledge-based setting,  there is no need for exhaustive and expensive (in time mainly) exploration of the state space, and the controller can directly choose training sets 
$\Gamma_n$ constructed from distributions $\mu_n$ that assign more points to the parts of the state space where the optimal process is likely to be driven. 
In practice, at time $n$, assuming we know that the optimal process is likely to stay in the ball centered around the point $m_n$ and with radius $r_n$, we choose a training measure $\mu_n$ centered around $m_n$ as, for example $\mathcal{N}(m_n,r_n^2)$, and build the training set as sample of the latter.

\vspace{3mm}

\noindent {\bf Explore first, exploit later}
If the agent has no idea of where to drive the process to receive large rewards, she can always proceed to an exploration step to discover favorable subsets of the state space.
To do so, $\Gamma_{n}$, the training sets at time $n$, for $n=0,\ldots,N-1$, can be built by forward simulation using optimal strategies estimated using algorithms developed e.g. in \cite{hanE16}. We are then back to the Exploitation only strategy framework.
An alternative suggestion is to use the approximate policy obtained from a first iteration of the hybrid-now algorithm, in order to derive the 
associated approximated optimal state process, $\hat X^{(1)}$. In a second iteration of the  hybrid-now algorithm, we then use  the marginal distribution $\mu_n^{(1)}$ of $\hat X_n^{(1)}$. 
We can even iterate this procedure in order to improve the choice of the  training measure $\mu_n$, as proposed  in \cite{kouetal16}.  

\vspace{3mm}

\begin{remark}
	\label{rmk:merge_trainings}
For stationary control problems, which are the common framework in reinforcement learning, it is usual to combine exploration of the (state and action) space and learning of the value function by drawing the mini-batches under the law given by the current estimate of the $Q$-value plus some noise (commonly called $\eps$-greedy policy). In practice, this method works well, and offers a flexible trade-off between exploration and exploitation. It looks challenging, however, to study its performance theoretically because it mixes the gradient-descent steps for the training of the value function and the optimal control. Also, note that this method is not easily adaptable to finite-horizon control problems. 
\end{remark}

\section{Convergence analysis} \label{secconv}

This section is devoted to the convergence of the estimator $\hat V_n^M$  of the value function $V_n$  obtained from a training sample of size $M$ and using DNN algorithms listed in Section \ref{sec::listeAlgos}.

Training samples rely on a given family of probability distributions $\mu_n$ on $\Xc$, for $n=0,\ldots,N$, referred to as training distribution (see Section \ref{sectraining} for a discussion on the choice of $\mu$). For the sake of simplicity, we consider that $\mu_n$ does not depend on $n$, and denote then by $\mu$ the training distribution.
We shall assume that the family of controlled transition probabilities has a density w.r.t. $\mu$, i.e., 
\beqs
P^a(x,dx') &=&  r(x,a;x') \mu(dx'). 
\enqs
We shall assume that $r$ is uniformly bounded in  $(x,x',a)$ $\in$ $\Xc^2\times\A$, and uniformly Lipschitz w.r.t. $(x,a)$, i.e., 

\vspace{1.5mm}

\noindent \setword{{\bf (Hd)}}{Hd} There exists some positive constants $\| r \|_\infty$ and  $[r]_{_L}$ s.t. 
\beqs
| r(x,a;x') | & \leq & \| r \|_\infty,  \;\;\; \forall x,x' \in \Xc, \; a \in \A, \\
| r(x_1,a_1;x') - r(x_2,a_2;x') | & \leq & [r]_{_L} ( |x_1-x_2| + |a_1 - a_2|), \;\;\; \forall x_1,x_2 \in \Xc, \; a_1,a_2 \in \A.
\enqs

\begin{remark}
\rm Assumption \ref{Hd} is usually satisfied when the state and control space are compacts. While the compactness on the control space $\A$ is not explicitly assumed, the compactness condition on the state space  $\Xc$ turns out to be more crucial for deriving estimates on the estimation error (see Lemma \ref{speedLearning}), and will be assumed to hold true for simplicity. Actually, this compactness condition on $\Xc$  can be relaxed by truncation and localization arguments (see proposition \ref{prop:localize} in Appendix A) by considering a training 
distribution $\mu$ such that \ref{Hd} is true and which admits a moment of order 1, i.e. $\int \vert y \vert d\mu(y)$ $<$ $+\infty$.
\ep
\end{remark}

\vspace{2mm}

\noindent We shall also assume some boundedness and Lipschitz condition on the reward functions:

\vspace{1.5mm}

\noindent\setword{{\bf (HR)}}{HR} There exists some positive constants $\| f \|_\infty$, $\| g \|_\infty$,  $[f]_{_L}$, and $[g]_{_L}$ s.t. 
\beqs
| f(x,a) | & \leq & \| f \|_\infty, \;\;\;  | g(x) | \; \leq \;  \| g \|_\infty,  \;\;\; \forall x \in \Xc, \; a \in \A, \\
| f(x_1,a_1) - f(x_2,a_2) | & \leq & [f]_{_L} ( |x_1-x_2| + |a_1 - a_2|), \\
| g(x_1) - g(x_2) | & \leq &  [g]_{_L}  |x_1-x_2|, \quad \;\;\; \forall x_1,x_2 \in \Xc, \; a_1,a_2 \in \A.
\enqs

\vspace{1mm}

Under this boundedness condition, it is clear that the value function $V_n$ is also bounded, and we have:
$
\lVert V_n \rVert_\infty \leq (N-n)\lVert f \rVert_\infty + \lVert g \rVert_\infty,$ for $ n \in \{0,...,N\}.
$

\vspace{2mm}

We shall finally assume a Lipschitz condition on the dynamics of the MDP. 

\vspace{1mm}

\noindent \setword{{\bf (HF)}}{HF} For any $e$ $\in$ $E$, there exists $C(e)$ such that for all couples $(x,a)$ and $(x',a')$ in $\Xc\times\A$:
\begin{equation}
\left\vert F(x,a,e) - F(x', a',e) \right\vert \leq C(e) \left( \vert x-x'\vert + \vert a-a' \vert   \right). \nonumber
\end{equation}

For  $M$ $\in$ $\N^*$, let us define $\rho_M = \E \left[ \Sup_{1 \leq m \leq M}C(\eps^m) \right], $
where the  $(\eps^m)_m$ is an i.i.d. sample of the noise $\eps$.  The rate of convergence of $\rho_M$ toward infinity will play a crucial role to show the convergence of the algorithms.

\begin{remark}
 {\rm  A typical example when  \ref{HF} holds  is the case where $F$ is defined through the time discretization of an Euler scheme, i.e. 
$
F(x,a,\eps) := b(x,a) + \sigma(x,a) \eps, \nonumber
$
with $b$ and $\sigma$ Lipschitz-continuous w.r.t. the couple $(x,a)$, and $\eps \sim \mathcal{N}(0,I_d)$, where $I_d$ is the identity matrix of size $d\times d$. 
Indeed, in this case, it is straightforward to see that $C(\eps)= [b]_L + [\sigma]_L \lVert \eps \rVert_d$, where $[b]_L$ and $[\sigma]_L$ stand for the Lipschitz coefficients of $b$ and $\sigma$, and $\lVert. \rVert_d$ stands for the Euclidean norm in $\R^d$. Moreover, one can show that:
\begin{equation}
\label{ineq:normalRV}
%\E \left[ \Sup_{1 \leq m \leq M} \prod_{k=n}^{N}C(\eps_k^m)\right] \leq
 \rho_M   \leq  [b]_L + d[\sigma]_L   \sqrt{2 \log(2dM)},
 \end{equation}
 which implies in particular that $
\rho_M  \underset{M \to +\infty}{=} \mathcal{O} \left(   \sqrt{\log(M)}  \right). \nonumber
$
Let us indeed check the inequality  \eqref{ineq:normalRV}. For this, let us fix some  integer $M'>0$ and let 
$Z:=\Sup_{1 \leq m \leq M'} \vert \epsilon_1^m \vert$ where $\epsilon_1^m$ are i.i.d. such that $\epsilon_1^1 \sim \mathcal{N}(0,1)$. 
From Jentzen inequality to the r.v. $Z$ and the convex function $z \mapsto \exp(t z)$, where $t>0$ will be fixed later, we get
{\small \begin{align}
e^{ t \E \left[ Z \right]} &\leq  \E \left[	\exp\big( t Z\big) \right] \leq \E \left[\Sup_{1 \leq m \leq M'} \exp \big(t \vert \epsilon_1^m \vert \big)  \right] 
\leq \sum_{m=1}^{M'}\E \Big[ \exp \big(t \vert \epsilon_1^m \vert \big)  \Big] \leq 2M'  \exp \big( \frac{t^2}{2} \big), \nonumber
\end{align}}
where we used the closed-form expression of the moment generating function of the folded normal distribution\footnote{The folded normal distribution is the one of $\vert Z \vert$ with $Z \sim \mathcal{N}_1(\mu,\sigma)$. Its moment generating function is given by $t \mapsto \exp\left(\frac{\sigma^2t^2}{2} + \mu t\right) \left[ 1- \Phi \left(-\frac{\mu}{\sigma} - \sigma t\right) \right] + \exp\left(\frac{\sigma^2t^2}{2} - \mu t\right) \left[ 1- \Phi \left(\frac{\mu}{\sigma} - \sigma t\right) \right]$, where $\Phi$ is the c.d.f. of $\mathcal{N}_1(0,1)$. } to write the last inequality.
Hence, we have for all $t>0$:
\begin{align}
\E\big[ Z\big] \leq \frac{\log(2M')}{t} + \frac{t}{2}. \nonumber
\end{align}
We get, after taking $t= \sqrt{2 \log(2M')}$:
\begin{align}
\E\big[ Z\big] \leq \sqrt{2 \log(2M')}. \label{ineq:majoreNormal}
\end{align}
Since inequality $\lVert x \rVert_d \leq d \lVert x \rVert_{\infty}$ holds for all $x \in \R^d$, we derive
$$
\E\left[ \Sup_{1 \leq m \leq M} C(\eps^m) \right] \leq [b]_L + d [\sigma]_L \E \left[ \Sup_{1 \leq m \leq dM} C(\epsilon_1^m) \right],
$$
and apply \eqref{ineq:majoreNormal} with $M'=dM$, to complete the proof of \eqref{ineq:normalRV}.

\ep
}
\end{remark}

\vspace{3mm}

\begin{remark}
{\rm 	
Under \ref{Hd},  \ref{HR} and \ref{HF}, it is straightforward to see from the dynamic programming formula \eqref{DP} 
that $V_n$ is Lipschitz for all $n=0,\ldots, N$, with a Lipschitz coefficient $[V_{n}]_L$ that can be bounded by standard arguments as follows: 
{\small	\begin{equation*}
	\left\{
	\begin{array}{ccl}
	[V_N]_L & = & [g]_L \\
	\left[ V_n \right]_L & \leq& \min \left(  [f]_L + \lVert V_n \rVert_\infty [r]_L \;, \;\rho_1 \frac{1- \rho_1^{N-n}}{1- \rho_1} + \rho_1^{N-n} [g]_L \right) , \text{ for } n=0,\ldots, N-1,
	\end{array}
	\right.
	\end{equation*}}
	using the usual convention $\frac{1-x^p}{1-x}=p$ for $p \in \mathbb{N}^*$ and $x=0$. 
	The Lipschitz continuity of $V_n$ plays a significant role in proving the convergence of the Hybrid algorithms described and studied in section \ref{subsectionHybrid}.
}
\ep	
\end{remark}

\subsection{Control learning  by performance iteration (NNcontPI)}	
%\label{subsectionLC}
\label{sec:NNContPI}

In this paragraph, we analyze the convergence of the NN  control learning  by performance iteration as described in Section \ref{secAlgoPI}. 
We shall consider the following set of neural networks for approximating the optimal policy: 
\beq
\quad \quad \quad ~^\eta \! \! \Ac_K^\gamma & := & \Big\{  x \in \Xc \mapsto A(x;\beta) \; = \; (A_1(x;\beta),\ldots,A_q(x;\beta)) \in \A, \; \label{defAetagamma} \\
& & \;\;\; A_i(x;\beta) = \sigma_{\A} \Big( \sum_{j=1}^{K} c_{ij} (a_{ij}.x+b_{ij})_+ + c_{0j} \Big), \;\;\;  i=1,\ldots,q,   \nonumber \\
& & \;\;\;\;\;\;\;  \beta = (a_{ij},b_{ij},c_{ij})_{i,j}, \; a_{ij} \in \R^d, \lVert a_{ij} \rVert \leq \eta,   b_{ij}, c_{ij} \in \R,  \sum_{i=0}^{K} \vert c_{ij} \vert \leq \gamma \Big\}, \nonumber 
\enq
where $\lVert . \rVert$ is the Euclidean norm in $\R^d$. Note that the considered neural networks have one hidden layer, $K$ neurons, Relu activation  functions, and $\sigma_{\A}$ as output layer. $\gamma$ and $\eta$ are often referred to in the literature as respectively total variation and kernel. 
The activation function $\sigma_A$ is chosen depending on the form of the control space $\mathbb{A}$. For example, when $\A$ $=$ $\R^q$, we simply take $\sigma_A$ as the identity function; when $\A$ 
$=$ $\R_+^q$, or more generally in the form  $\A$ $=$ $\Prod_{i=1}^q [a_i,\infty)$,  
one can take  the component-wise  Relu activation function (possibly shifted and scaled); and when $\A$ $=$ $[0,1]^q$, or more generally 
$\A$ $=$  $\Prod_{i=1}^q [a_i,b_i]$, for $a_i$ $\leq$ $b_i$, $i$ $=$ $1,\ldots,q$, one can take  the component-wise  sigmoid  activation function (possibly shifted and scaled); 
when $\A$ is a discrete space corresponding to discrete control values $\{a_1,\cdots,a_L\}$, we may consider randomized control, hence valued in the simplex of $\hat\A$ of  $[0,1]^L$, 
and use a softmax output layer activation function for parametrizing the discrete probabilities valued in $\hat\A$ (such algorithm, called ClassifPI, has been developed in \cite{bacetal18}; 
see also \cite{ger17} on softmax output for classification problems, and \cite{asadi17} about softmax policies for action selection and possible improvements.)

\vspace{3mm}

Let $K_M$, $\eta_M$ and $\gamma_M$ be sequences of integers such that 	
\begin{equation}
\begin{split}
K_M, \gamma_M, \eta_M \xrightarrow[M \to \infty]{} \infty \quad \text{and} \quad \rho_M^{N-1}\gamma_M^{N-1} \eta_M^{N-2}\sqrt{\frac{\log(M)}{M}} \xrightarrow[M \to \infty]{} 0.
\end{split}
\label{cond::algoIII}
\end{equation}
We denote by  $\mathcal{A}_M$ $:=$ $~^{\eta_M}\! \!\mathcal{A}_{K_M}^{\gamma_M}$ the class of neural network  for policy as defined in \eqref{defAetagamma} 
with parameters $\eta_M$, $K_M$ and $\gamma_M$ that satisfy the conditions \eqref{cond::algoIII}. 
\begin{remark}
{\rm 
In the case where $F$ is defined in dimension $d$ as:
$
F(x,a,\eps)= b(x,a) + \sigma(x,a) \eps,
$
we can use \eqref{ineq:normalRV} to get: 
$
\rho_M^{N-n} \underset{M \to +\infty }{=} \mathcal{O}\left( \sqrt{\log(M)}^{N-n}\right).
$
\ep
}
\end{remark}

Recall that the  approximation of the optimal policy in the NNcontPI algorithm is computed in backward induction as follows: 
For $n$ $=$ $N-1,\ldots,0$,  generate a training sample for the state $X_n^{(m)}$, $m$ $=$ $1,\ldots,M$  from the training distribution $\mu$, and 
samples of the exogenous noise $\big( \eps_k^m \big)_{m=1,k=n+1}^{M,N}$.
 	\begin{itemize}
		\item Compute the approximated policy at time $n$
		\begin{equation}
		\label{eq:PI_optim}
		\begin{array}{ccc}
		\hat{a}_n^M & \in &  \argmin_{A \in \mathcal{A}_M} \frac{1}{M} \sum_{m=1}^{M} \left[ f(X_n^{(m)},A(X_n^{(m)}))+ \hat Y^{(m),A}_{n+1}\right]
		\end{array}
		\end{equation}
		where 
		\begin{align} \label{defY}
			\hat Y^{(m),A}_{n+1}= \sum_{k=n+1}^{N-1} f\left( X^{(m),A}_k, \hat{a}_k^M\left( X^{(m),A}_k \right)\right)  + g\left(  X^{(m),A}_N \right),
		\end{align}
		with $ \big(X^{(m),A}_k\big)_{k=n+1}^N$ defined by induction as follows, for $m$$=$$1,\ldots,M$:
		\begin{equation*}
		\left\{
		\begin{array}{ccl}
		X^{(m),A}_{n+1} &=& F\Big( X^{m}_n, A\big(X^{m}_n \big), \eps_{n+1}^m \Big) \\
	        X^{(m),A}_{k} &=& F\Big( X^{(m),A}_{k-1}, \hat{a}_k^M\big(X^{(m),A}_{k-1} \big), \eps_{k}^m \Big), \;\;\; \mbox{  for } k= n+2, \ldots, N.
	        \end{array}
	        \right.
		\end{equation*}
%\vspace{3mm}
		\item Compute the estimated value function $\hat V_n^M$ as in \reff{estimValuePI}. 
	\end{itemize}

%\bl{parler des résultats de Noufel sur la descente de gradient sto et l'erreur associée.}
\begin{remark} 
{\rm Finding the $\argmin_{}$ in \eqref{eq:PI_optim} is highly non-trivial and impossible in real situation. It is slowly approximated by (non stochastic) gradient-descent, i.e. with one large batch of size $M$, for which some convergence results are known under convexity or smoothness assumptions of $F$, $f$, $g$; and is approximated in practice by stochastic gradient-descent with a large enough number of epochs. We refer to \cite{jenetal18} for error analysis of stochastic gradient descent optimization algorithms. 
	In order to simplify the theoretical analysis, we assume that the $\argmin_{}$ in \eqref{eq:PI_optim} is exactly reached. Our main focus here is to understand how the finite size of the training set affects the convergence of our algorithms.
}
\ep	
\end{remark}

\vspace{2mm}

We now state our main result about the convergence of the NNcontPI algorithm.

\begin{theorem}	\label{theo::algoIII}
Assume that there exists an optimal feedback control $(a_k^{\mathrm{opt}})_{k=n}^{N-1}$ for the control problem with value function $V_n$, $n$ $=$ $0,\ldots,N$, and  
let $X_n \sim \mu$. Then, as $M \to \infty$\footnote{The notation $x_M$ $=$ $\Oc(y_M)$ as $M$ $\rightarrow$ $\infty$, means that the ratio $|x_M|/|y_M|$ is bounded as $M$ goes to infinity.} 
\begin{equation}
\label{ineq::restheoIII}
\begin{split}
\E \big[  \hat{V}_n^M(X_n) - V_n(X_n) \big]  & =  
\mathcal{O} \Bigg(\frac{\rho_M^{N-n-1}\gamma_M^{N-n-1} \eta_M^{N-n-2}}{\sqrt{M}} \\
 & \hspace{.9cm} + \Sup_{n \leq k \leq N-1}\inf_{A \in \mathcal{A}_M} \E \left[ \vert A(X_k) - a^{\mathrm{opt}}_k(X_k) \vert  \right]  \Bigg),
\end{split}
\end{equation}
where $\E$ stands for the expectation over the training set used to evaluate the approximated optimal policies $(\hat a_k^M)_{n \leq k \leq N-1}$, 
as well as the path $(X_n)_{n \leq k\leq N}$ controlled by the latter. 
Moreover, as $M \to \infty$\footnote{The notation $x_M$ $=$ $\Op (y_M)$ as $M$ $\rightarrow$ $\infty$, means that there exists $c>0$ such that $\P \big(\vert x_M \vert >c \vert y_M \vert \big) \to 0 $ as $M$ goes to infinity.} 
\begin{equation}
\begin{split}
	\E_M \big[  \hat{V}_n^M(X_n) - V_n(X_n) \big]  &= 
\Op \Bigg(\rho_M^{N-n-1}\gamma_M^{N-n-1} \eta_M^{N-n-2}\sqrt{\frac{\log(M)}{M}} \\
& \hspace{1.2cm}+ \sup_{n \leq k \leq N-1}\inf_{A \in \mathcal{A}_M} \E \left[ \vert A(X_k) - a^{\mathrm{opt}}_k(X_k) \vert  \right]  \Bigg), \label{ineq::restheoIII_2}
\end{split}
\end{equation}
where $\E_M$ stands for the expectation conditioned by the training set used to estimate the optimal policies $(\hat a_k^M)_{n \leq k \leq N-1}$.
\end{theorem}

\begin{remark}
{\rm
{\bf 1.}  The term $\frac{\rho_M^{N-n-1}\gamma_M^{N-n-1} \eta_M^{N-n-2}}{\sqrt{M}}$  should be seen as the estimation error. It is due to the approximation of  the optimal controls by means of neural networks in $\mathcal{A}_M$ using \emph{empirical} cost functional in \reff{eq:PI_optim}. We show  in section \ref{sec::forwardEvaluation} that this term disappears in the ideal case where the real cost functional  (i.e. not the empirical one) is minimized.

{\bf 2.} The rate of convergence depends dramatically on $N$ since it becomes exponentially slower when $N$ goes to infinity. This is a huge drawback for this performance iteration-based algorithm. We will see in the next section that the rate of convergence of value iteration-based algorithms does not suffer from this problem.
}
\ep
\end{remark}

\noindent {\bf Comment}: 
Since we clearly have  $V_n$ $\leq$ $\hat V_n^M$, estimation 
\eqref{ineq::restheoIII} implies the  convergence  in $L^1$ norm of the NNcontPI algorithm, under condition \eqref{cond::algoIII}, and in the case where 
$${\sup_{n \leq k \leq N}\inf_{A \in \mathcal{A}_M} \E \big[ \vert A(X_k) - a^{\mathrm{opt}}_k(X_k) \vert  \big] \xrightarrow[M \to +\infty]{} 0}\ .$$ 
This is actually the case under some regularity assumptions on the optimal controls, as stated in the following proposition.

\begin{proposition}	\label{ConvergenceError}
The two following assertions hold:
	\begin{enumerate}
		\item Assume that $a^{\mathrm{opt}}_k \in  \mathbb{L}^1(\mu)$ for $k=n,...,N-1$. Then
		\begin{equation}
		\label{eq:denseness1}
		\sup_{n \leq k \leq N-1}\inf_{A \in \mathcal{A}_M} \E \big[ \vert A(X_k) - a^{\mathrm{opt}}_k(X_k) \vert  \big] \xrightarrow[M \to +\infty]{} 0.
		\end{equation}
		\item  Assume that the function $a^{\mathrm{opt}}_k$ is  $c$-Lipschitz for $k=n,...,N-1$. Then 
		\begin{align}
		\sup_{n \leq k \leq N-1}\inf_{A \in \mathcal{A}_M} \; \E \big[ \vert A(X_k) - a^{\mathrm{opt}}_k(X_k) \vert  \big]  &  \nonumber\\
		& \hspace{-3cm}< c \left( \frac{\gamma_M}{c}\right)^{-2d/(d+1)}\log \left(\frac{\gamma_M}{c} \right) + \gamma_M K_M^{-(d+3)/(2d)}.
		\label{eq:denseness2}
		\end{align}
	\end{enumerate}
\end{proposition}
{\bf Proof.}  The first statement of Proposition \ref{ConvergenceError} relies essentially on the universal approximation theorem, 
and the second assertion is stated and proved in \cite{bach14}. For the sake of completeness, we provide more details in  Appendix \ref{sec::functionApproximationUsingNeuralNetworks}.
\ep

\begin{remark}
The consistency of NNContPI is proved in Theorem \ref{theo::algoIII}, but no results on its variance are provided. We refer e.g. to Section 3.2 in \cite{bacetal18} to observe the small variance of NNContPI for a linear quadratic problem in dimension up to 100.
\end{remark}

\vspace{1mm}

The rest of this section is devoted to the proof of Theorem \ref{theo::algoIII}.  Let us introduce some useful notations. Denote by $\A^{\Xc}$ the set of Borelian functions from the state space $\Xc$ into the control space $\A$.  For $n$ $=$ $0,\ldots,N-1$, and given a feedback control (policy) represented by  a sequence $(A_k)_{k=n,\ldots,N-1}$, with $A_k$ in  $\A^{\Xc}$, we denote by $J_n^{(A_k)_{k=n}^{N-1}}$ the cost functional associated with the policy $(A_k)_k$.  Notice that with this notation, we have $\hat V_n^M$ $=$ 
$J_n^{(\hat a_k^M)_{k=n}^{N-1}}$. We define the {\it estimation error} at time $n$ associated with the NNContPI algorithm by
\beqs 
\eps_{\mathrm{PI},n}^{\mathrm{esti}} & := &  \sup_{A \in \Ac_M}  \Big|  \frac{1}{M} \sum_{m=1}^M \Big[ f(X_n^{(m)},A(X_n^{(m)})) +  \hat Y_{n+1}^{(m),A} \Big]  - 
\mathbb{E}_M \big[ J_n^{A,(\hat a_k^M)_{k=n+1}^{N-1}}(X_n) \big] \Big|, 
\enqs
with $X_n$ $\sim$ $\mu$: It measures how well the chosen estimator (e.g. mean square estimate) can approximate a certain quantity (e.g. the conditional expectation). Of course we expect the latter to cancel when the size of the training set used to build the estimator goes to infinity.  Actually, we have

\begin{lemma} \label{speedLearning} For $n$ $=$ $0,\ldots,N-1$, the following holds:
{\small \begin{align}  & \hspace{-5mm}\E[\eps_{\mathrm{PI},n}^{\mathrm{esti}}]\leq\big(\sqrt{2}+16\big)\frac{\big((N-n)\lVert f\rVert_{\infty}+\lVert g\rVert_{\infty}\big)}{\sqrt{M}} \nonumber\\
&\!+\!\frac{16\gamma_{M}}{\sqrt{M}}\Bigg\{[f]_{L}\!\left(1\!+\!\rho_{M}\frac{1\!-\!\rho_{M}^{N-n-1}\big(1\!+\!\eta_{M}\gamma_{M}\big)^{N-n-1}}{1\!-\rho_{M}\big(1\!+\!\eta_{M}\gamma_{M}\big)}\right)\!+\big(1\!+\!\eta_{M}\gamma_{M}\big)^{N-n-1}\rho_{M}^{N-n}[g]_{L}\Bigg\}\nonumber\\
& \hspace{5mm} =\Oc\bigg(\frac{\rho_{M}^{N-n-1}\gamma_{M}^{N-n-1}\eta_{M}^{N-n-2}}{\sqrt{M}}\bigg),\quad\text{ as \ensuremath{M\to\infty}}.\label{estimesti}
\end{align}}
This implies in particular that 	
\begin{align}
 \label{estimesti_2}
\eps_{\mathrm{PI},n}^{\mathrm{esti}}  & = \Op  \left(\rho_M^{N-n-1}\gamma_M^{N-n-1} \eta_M^{N-n-2} \sqrt{\frac{\log(M)}{M}}\right), \quad \text{ as $M \to \infty$},
\end{align}
where we recall that $\rho_M = \E\left[ \Sup_{1 \leq m \leq M} C(\eps^m ) \right]$ is defined in \ref{HF}.
\end{lemma}
{\bf Proof.}
The relation \eqref{estimesti} states that the estimation error cancels when $M\to \infty$ with a rate of convergence of order  $\mathcal{O}\Big (\frac{\rho_M^{N-n-1}\gamma_M^{N-n-1} \eta_M^{N-n-2}}{\sqrt{M}} \Big)$. The proof is in the spirit of the one that can be found in chapter 9 of \cite{gyokoh02}. 
 It relies on a technique of symmetrization by a ghost sample, and a wise introduction of additional randomness by random signs. 
 The details are postponed to Section \ref{sec::estimationError} in the Appendix.  
 The proof of \reff{estimesti_2} follows from \reff{estimesti} by a direct application of the Markov inequality. 
\ep

\vspace{3mm}

Let us also define the {\it approximation error} at time $n$ associated with the NNContPI algorithm by
 \begin{align}\label{defepsapprox}
 \eps_{\mathrm{PI},n}^{\mathrm{approx}} & := 
 \inf_{A \in \mathcal{A}_M} \mathbb{E}_M \Big[ J_n^{A,(\hat a_k^M)_{k=n+1}^{N-1}}(X_n) \Big] 
 - \inf_{A \in \A^\Xc} \mathbb{E}_M \Big[ J_n^{A,(\hat a_k^M)_{k=n+1}^{N-1}}(X_n) \Big],
 \end{align}
 where we recall that $\E_M$ denotes the expectation conditioned by the training set used to compute the estimates $(\hat a_k^M)_{k=n+1}^{N-1}$ and the one of $X_n$ $\sim$ $\mu$.\\
 $\eps_{\mathrm{PI},n}^{\mathrm{approx}}$ measures how well the regression function can be approximated by means of neural networks functions in $\mathcal{A}_M$ 
 (notice that the class of neural networks is not dense in the set $\A^{\Xc}$ of all Borelian functions).

 \begin{lemma} \label{lemapprox}
 	For $n$ $=$ $0,\ldots,N-1$, the following holds as $M \to \infty$,		
{ \small 	\begin{equation} 
 	\E [\eps_{\mathrm{PI},n}^{\mathrm{approx}}] = \mathcal{O} \left( \frac{\rho_M^{N-n-1}\gamma_M^{N-n-1} \eta_M^{N-n-2}}{\sqrt{M}}+ \sup_{n \leq k \leq N-1}\inf_{A \in \mathcal{A}_M} \E \big[ \vert A(X_k) - a^{\mathrm{opt}}_k(X_k) \vert  \big] \right).\label{ineq::BoundApproxError}
 	\end{equation}}
This implies in particular  	
{\begin{align} \label{ineq::BoundApproxError_2}
\eps_{\mathrm{PI},n}^{\mathrm{approx}}  & = \Op \left( {\rho_M^{N-n-1}\gamma_M^{N-n-1} \eta_M^{N-n-2}}\sqrt{\frac{\log(M)}{M}}  \right. \nonumber\\ 
& \hspace{1.2cm}+  \left.\sup_{n \leq k \leq N-1}\inf_{A \in \mathcal{A}_M} \E \big[ \vert A(X_k) - a^{\mathrm{opt}}_k(X_k) \vert  \big] \right).
\end{align}}
 \end{lemma}
 {\bf Proof.} See  Appendix \ref{sec::annexerror1}  for the proof of \eqref{ineq::BoundApproxError}. \eqref{ineq::BoundApproxError_2} is  a direct application of the Markov inequality.  	 
 \ep  

\vspace{3mm}

\noindent {\bf Proof  of Theorem \ref{theo::algoIII}.} \\
 {\it Step 1.} Let us denote by 
$
\hat{J}_{n,M}^{A,(\hat{a}^M_{k})_{k=n+1}^{N-1}} :=  \frac{1}{M} \sum_{m=1}^{M} \Big[  f\big(X_n^{(m)},A(X_n^{(m)}) \big)+\hat Y_{n+1}^{(m),A}  \Big], 
$
 the empirical cost function, from time $n$ to $N$,  associated with the sequence of controls $(A,(\hat{a}^M_k)_{k=n+1}^{N-1})$ and the training set, where we recall that $\hat Y_{n+1}^{(m),A}$ is defined in \reff{defY}.  
 We then have 
 \begin{align}
 \label{334}
 \E_M\big[\hat{V}^M_n(X_n)\big] &= \E_M\Big[J_n^{(\hat a_k^M)_{k=n}^{N-1}}(X_n)\Big]- \hat{J}_{n,M}^{(\hat{a}^M_k)_{k=n}^{N-1}} + \hat{J}_{n,M}^{(\hat{a}^M_k)_{k=n}^{N-1}}  \nonumber\\
% &\leq   \sup_{A \in \mathcal{A}_M} \bigg\vert \E_M\Big[J_n^{A,(\hat a_k^M)_{k=n+1}^{N-1}}(X_n)\Big]- \hat{J}_{n,M}^{A,(\hat{a}^M_{k})_{k=n+1}^{N-1}} 
 %\bigg\vert +\hat{J}_{n,M}^{(\hat{a}^M_k)_{k=n}^{N-1}} \nonumber \\
 &\leq  \eps_{\mathrm{PI},n}^{\mathrm{esti}} +\hat{J}_{n,M}^{(\hat{a}^M_k)_{k=n}^{N-1}}, 
 \end{align}
 by definition of $\hat V_n^M$ and $\eps_{\mathrm{PI},n}^{\mathrm{esti}}$. Moreover, for any $A \in \mathcal{A}_M$,
 \begin{align}
  \label{ineq::intermediaire}
\hat{J}_{n,M}^{A,(\hat{a}^M_{k})_{k=n+1}^{N-1}} &= 
\hat{J}_{n,M}^{A,(\hat{a}^M_{k})_{k=n+1}^{N-1}} -  \E_M\Big[J_n^{A,(\hat a_k^M)_{k=n+1}^{N-1}}(X_n)\Big]+ \E_M\Big[J_n^{A,(\hat a_k^M)_{k=n+1}^{N-1}}(X_n)\Big] \nonumber\\
&\leq    \eps_{\mathrm{PI},n}^{\mathrm{esti}} + \E_M\Big[J_n^{A,(\hat a_k^M)_{k=n+1}^{N-1}}(X_n)\Big]. 
\end{align}
Recalling that 
$
\hat{a}_n^M  =   \argmin_{A \in \mathcal{A}_M} \hat{J}_{n,M}^{A,(\hat{a}^M_{k})_{k=n+1}^{N-1}},
$
and  taking   the infimum over $\mathcal{A}_M$ in the l.h.s. of \eqref{ineq::intermediaire} first, and  in the  r.h.s. second, we then get
 \beqs 
 \hat{J}_{n,M}^{(\hat{a}^M_k)_{k=n}^{N-1}} &\leq &  \eps_{\mathrm{PI},n}^{\mathrm{esti}} + \inf_{A \in \mathcal{A}_M} \E_M\Big[J_n^{A,(\hat a_k^M)_{k=n+1}^{N-1}}(X_n)\Big].
 \enqs
Plugging this last inequality into \eqref{334}  yields the following estimate
\begin{equation} \label{interestistep1}
\begin{array}{ccc}
\E_M\big[\hat{V}^M_n(X_n)\big]   -  \inf_{A \in \mathcal{A}_M} \E_M\Big[ J_n^{A,(\hat a_k^M)_{k=n+1}^{N-1}}(X_n) \Big] &\leq & 2\eps_{\mathrm{PI},n}^{\mathrm{esti}}.
\end{array}
\end{equation}

\noindent {\it Step 2.}   By definition \reff{defepsapprox} of the approximation error, using the law of iterated conditional expectations for $J_n$, and the dynamic programming principle for $V_n$ with the optimal control $a_n^{\mathrm{opt}}$ at time $n$, we have
\begin{align*}
\inf_{A \in \mathcal{A}_M}  \E_M\big[ J_n^{A,(\hat a_k^M)_{k=n+1}^{N-1}}(X_n) \big] -  \E_M[ V_n(X_n) ]   \\
& \hspace{-4cm}=   \eps_{\mathrm{PI},n}^{\mathrm{approx}} + \inf_{A \in \A^{\Xc}} \E_M \big\{ f(X_n,A(X_n)) + \E_n^{A}\big[ J_{n+1}^{(\hat{a}_k^M)_{k=n+1}^{N-1}} (X_{n+1}) \big] \big\}  \\
&  \hspace{-3.4cm}-  \E_M \Big[ f(X_n,a^{\mathrm{opt}}_n(X_n))  +  \E_n^{a^{\mathrm{opt}}_n}\big[V_{n+1}(X_{n+1})\big]  \Big]    \\
& \hspace{-4cm} \leq   \eps_{\mathrm{PI},n}^{\mathrm{approx}}  \; + \;    \E_M  \E_n^{a^{\mathrm{opt}}_n} \Big[ J_{n+1}^{(\hat{a}_k^M)_{k=n+1}^{N-1}}(X_{n+1})-V_{n+1}(X_{n+1}) \Big], 
\end{align*}
where $\E_n^{A}[.]$ stands for the expectation conditioned by $X_n$ at time $n$ and the training set, when strategy $A$ is followed at time $n$. Under the bounded density assumption in \ref{Hd}, we then get
\begin{align}\label{interestistep2} 
 \inf_{A \in \mathcal{A}_M}  \E_M\big[ J_n^{A,(\hat a_k^M)_{k=n+1}^{N-1}}(X_n) \big] -  \E_M[ V_n(X_n) ]& \nonumber \\
& \hspace{-6cm} \leq    \eps_{\mathrm{PI},n}^{\mathrm{approx}}  \; + \; \| r \|_{\infty} \int  \big[  J_{n+1}^{(\hat{a}_k^M)_{k=n+1}^{N-1}}(x')-V_{n+1}(x')  \big] \mu(dx') \nonumber \\
&  \hspace{-6cm} \leq    \eps_{\mathrm{PI},n}^{\mathrm{approx}}  \; + \; \| r \|_{\infty} \mathbb{E}_M\Big[ \hat{V}_{n+1}^M(X_{n+1}) - V_{n+1}(X_{n+1}) \Big], \text{ with } X_{n+1} \sim \mu.
\end{align}

\vspace{1mm}

\noindent {\it Step 3.} From \eqref{interestistep1} and \eqref{interestistep2}, we have 
\begin{align}
\label{interstep3}
 & \mathbb{E}_{M}\big[\hat{V}_{n}^{M}(X_{n})-V_{n}(X_{n})\big]=\E_{M}\big[\hat{V}_{n}^{M}(X_{n})\big]-\inf_{A\in\mathcal{A}_{M}}\E_{M}\big[J_{n}^{A,(\hat{a}_{k}^{M})_{k=n+1}^{N-1}}(X_{n})\big]\Big]\nonumber\\
 & +\inf_{A\in\mathcal{A}_{M}}\E_{M}\big[J_{n}^{A,(\hat{a}_{k}^{M})_{k=n+1}^{N-1}}(X_{n})\big]-\E_{M}[V_{n}(X_{n})]\nonumber\\
 & \leq2\eps_{\mathrm{PI},n}^{\mathrm{esti}}+\eps_{\mathrm{PI},n}^{\mathrm{approx}}+\|r\|_{\infty}\mathbb{E}_{M}\Big[\hat{V}_{n+1}^{M}(X_{n+1})-V_{n+1}(X_{n+1})\Big].
\end{align}
By induction, this implies 
\begin{align*}
%\label{interestistep2bis}
\mathbb{E}_M\big[ \hat{V}_n^M(X_n) - V_n(X_n) \big] & \leq  \sum_{k=n}^{N-1} \big(   2\eps_{\mathrm{PI},k}^{\mathrm{esti}}  + \eps_{\mathrm{PI},k}^{\mathrm{approx}} \big).
\end{align*}
Use the estimations \eqref{estimesti_2} for $\eps_{\mathrm{PI},n}^{\mathrm{esti}}$ in Lemma \ref{speedLearning}, and \reff{ineq::BoundApproxError_2} for $\eps_{\mathrm{PI},n}^{\mathrm{approx}}$ in Lemma  
\ref{lemapprox}, and observe that  $\hat{V}_n(X_n)$ $\geq$ $V_n(X_n)$ holds a.s., to complete the proof of \reff{ineq::restheoIII_2}. 
Finally,  the proof of \reff{ineq::restheoIII} is obtained by taking expectation in \eqref{interstep3}, and using estimations \reff{estimesti} and \reff{ineq::BoundApproxError}.  
\ep

\subsection{Hybrid-Now algorithm}	\label{subsectionHybrid}

In this paragraph, we analyze the convergence of the hybrid-now  algorithm. 
We shall consider the following set of neural networks for the value function approximation:
\beqs
~^\eta \!\mathcal{V}_K^{\gamma} & := & \Big\{ x \in \Xc \mapsto \Phi(x;\theta) \; = \; \sum_{i=1}^{K} c_i \sigma(a_i.x+b_i)+c_0, \;   \\
& &  \;\;\;\;\;  \theta = (a_i,b_i,c_i)_i, \;\;  \;  \lVert a_i \rVert \leq \eta, \; b_i \in \mathbb{R}, \; \sum_{i=0}^{K} \vert c_i \vert \leq \gamma  \Big\}.
\enqs

Let $\eta_M$, $K_M$ and $\gamma_M$ be integers such that: 
\begin{equation}
	\eta_M, \gamma_M, K_M \xrightarrow[M \to \infty]{} \infty\; s.t. \;
	 \frac{\gamma_M^4 K_M \log(M)}{M}, \frac{ \gamma_M^4 \rho_M^2 \eta_M^2 \log(M)}{M}  \xrightarrow[M \to \infty]{} 0, \label{cond::hybrid}
\end{equation}
with $\rho_M$ defined in \ref{HF}.

In what follows we denote by $\mathcal{V}_M:=~^{\eta_M}  \! \! \mathcal{V}_{K_M}^{\gamma_M}$ the space of neural networks for the estimated value functions at time $n$ $=$ $0,\ldots,N-1$, parametrized by the values $\eta_M$, $\gamma_M$ and $K_M$ that satisfy \eqref{cond::hybrid}. We also consider the class  $\Ac_M$  of neural networks for  estimated feedback optimal control at time $n$ $=$ $0,\ldots,N-1$, as described in Section \ref{sec:NNContPI}, with the same parameters  $\eta_M$, $\gamma_M$ and $K_M$. 
\vspace{3mm}

Recall that the  approximation of the value function and optimal policy in the hybrid-now algorithm is computed in backward induction as follows: 
\begin{itemize}
\item Initialize $\hat V_N^M$ $=$ $g$
\item For $n$ $=$ $N-1,\ldots,0$,  generate a training sample $X_n^{(m)}$, $m$ $=$ $1,\ldots,M$ from the training distribution $\mu$, and a training sample for the exogenous noise $\eps_{n+1}^{(m)}$, $m$ $=$ $1,\ldots,M$.
\begin{itemize}
\item[(i)] compute the approximated policy at time $n$
\beqs
\hat{a}_n^M & \in &  \argmin_{A \in \mathcal{A}_M} \frac{1}{M} \sum_{m=1}^{M} \big[ f(X_n^{(m)},A(X_n^{(m)}))+ \hat{V}_{n+1}^M(X_{n+1}^{(m),A}) \big]
\enqs
where $X_{n+1}^{(m),A}$ $=$ $F(X_n^{(m)},A(X_n^{(m)}),\eps_{n+1}^{(m)})$ $\sim$ $P^{A(X_n^{(m)})}(X_n^{(m)},dx')$.  
\item[(ii)]  compute the untruncated estimation of the value function at time $n$
\beqs
\tilde V_n^M & \in & \argmin_{\Phi \in \mathcal{V}_M} \frac{1}{M} \sum_{m=1}^{M} \Big[  f(X_n^{(m)},\hat a_n^M(X_n^{(m)}))+ \hat{V}_{n+1}^M(X_{n+1}^{(m),\hat a_n^M}) - \Phi(X_n^{(m)}) \Big]^2
\enqs
and set the truncated estimated value function at time $n$
\beq \label{trunca}
\hat V_n^M & = &  \max \Big(  \min \big(  \tilde V_n^M , \lVert V_n \rVert_\infty \big) ,  -\lVert V_n \rVert_\infty \Big).
\enq
%where we use the notation $a \wedge b := \min(a,b)$ and $ a \vee b := \max(a,b)$. 
\end{itemize}
\end{itemize}

\vspace{2mm}	 

\begin{remark}
{\rm 	Notice that  we have truncated  the estimated value function in \reff{trunca} by an  {\it a priori} bound on the true value function.  
This truncation step is natural from a practical implementation point of view, and is also used for simplifying the proof of the convergence of the algorithm. 
}
\ep
\end{remark}

\vspace{2mm}

We now state our main result  about the convergence of the Hybrid-Now  algorithm.

\begin{theorem}\label{theoKohler}
Assume that there exists an optimal feedback control $(a_k^{\mathrm{opt}})_{k=n}^{N-1}$ for the control problem with value function $V_n$, $n$ $=$ $0,\ldots,N$, and  
let $X_n \sim \mu$. Then, as ${M\rightarrow+\infty}$
	\begin{align}
	\small
		\label{eq::theoHybrid}
\E_M\Big[\vert  \hat{V}_n^M(X_n)-V_n(X_n) \vert  \Big]& \nonumber \\
& \hspace{-2cm}= \mathcal{O}_\mathbb{P}\Bigg( \left( \gamma_M^4 \frac{K_M\log(M)}{M} \right)^{\frac{1}{2(N-n)}} + \left( \gamma_M^4 \frac{\rho_M^2 \eta_M^2 \log(M)}{M} \right)^{\frac{1}{4(N-n)}} \nonumber\\
& \hspace{-1.cm}+ \sup_{n \leq k \leq N} \inf_{\Phi \in \mathcal{V}_M} \left(\E_M\Big[ \vert \Phi(X_k)-V_k(X_k) \vert^2 \Big] \right)^{\frac{1}{2(N-n)}}  \nonumber \\
& \hspace{-1.cm}+\sup_{n \leq k \leq N}\inf_{A \in \mathcal{A}_M} \left(  \mathbb{E}\Big[ \vert A(X_k)-a^{\mathrm{opt}}_k(X_k) \vert \Big] \right)^{\frac{1}{2(N-n)}}\Bigg),
\end{align}
where $\E_M$ stands for the expectation conditioned by the training set used to estimate the optimal policies $(\hat a_k^M)_{n \leq k \leq N-1}$. 
\end{theorem}

\begin{remark}
The consistency of Hybrid-Now is proved in Theorem \ref{theoKohler}, but no results on its variance are provided. The stability of Hybrid-Now has nevertheless been observed numerically in  \cite{bacetal18}, and we refer e.g. to its Section 3.2 to observe the small variance of Hybrid-Now in a linear-quadratic problem in dimension up to 100.
\end{remark}

\noindent \textbf{Comment:} Theorem \ref{theoKohler} states that the estimator for the value function provided by the hybrid-now algorithm converges in 
$\mathbb{L}^1(\mu)$ when the size of the training set goes to infinity. 
Note that $\left( \gamma_M^4 \frac{K_M\log(M)}{M} \right)^{\frac{1}{2(N-n)}} $ and $\left( \gamma_M^4 \frac{\rho_M^2 \eta_M^2 \log(M)}{M} \right)^{\frac{1}{4(N-n)}}$ stand for the estimation error made by estimating \emph{empirically} respectively the value functions and the optimal control by neural networks; also: $\Sup_{n \leq k \leq N}\Inf_{\Phi \in \mathcal{V}_M} \sqrt{\mathbb{E}\Big[ \vert \Phi(X_k)-V_k(X_k) \vert^2 \Big]} $ and $\Sup_{n \leq k \leq N}\Inf_{A \in \mathcal{A}_M} \mathbb{E}\Big[ \vert A(X_k)-a^{\mathrm{opt}}_k(X_k) \vert \Big] $are the approximation error made by estimating  respectively the value function  and the optimal control by neural networks.

\vspace{3mm}

In order to prove Theorem \ref{theoKohler}, let us first introduce the estimation error at time $n$ associated with the Hybrid-Now algorithm by
\begin{align*}
\eps_{\mathrm{HN},n}^{\mathrm{esti}}  &:=  \sup_{A \in \Ac_M}  \Bigg|  \frac{1}{M} \sum_{m=1}^M \left[ f(X_n^{(m)},A(X_n^{(m)})) +  \hat Y_{n+1}^{(m),A} \right] \\
&\hspace{2cm}  - \mathbb{E}_{M,n,X_n}^A \left[f(X_n,A(X_n)) +   \hat{V}_{n+1}^M\big(X_{n+1}\big) \right] \Bigg|, 
\end{align*}
where $\hat Y_{n+1}^{(m),A} = \hat{V}_{n+1}^M\big(X_{n+1}^{(m),A}\big)$,
and $X_{n+1}^{(m),A}= F\left( X_n^{(m)},A(X_n^{(m)}),\eps_{n+1}^m\right)$.
We have the following bound on this estimation error:
 
 \vspace{1mm}

\begin{lemma} \label{speedLearning_Hybrid}
For $n$ $=$ $0,...,N-1$, the following holds:
	\begin{align} 
	\E \left[\eps_{\mathrm{HN},n}^{\mathrm{esti}} \right]  & \leq    \frac{ \big(  \sqrt{2} + 16 \big) \big( (N-n) \lVert f \rVert_{\infty} + \lVert g \rVert_{\infty} \big)+16 [f]_L}{\sqrt{M}} +16  \frac{\rho_M \eta_M \gamma_M^2}{\sqrt{M}} \nonumber \\
	&\underset{M \to \infty}{=} \Oc \bigg( \frac{\rho_M \eta_M \gamma_M^2  }{\sqrt{M}}\bigg).  	\label{estimesti_hybrid} 
	\end{align}
\end{lemma}
{\bf Proof.}
See Appendix \ref{sec::estimationError_VI}.
\ep

\begin{remark}
{\rm
The result stated by lemma \ref{speedLearning_Hybrid}  is  sharper   than the one stated in Lemma \ref{speedLearning} for the performance iteration procedure. The main reason is that we can make use of the  $\gamma_M \eta_M$-Lipschitz-continuity of the estimate of the value function at time 
$n+1$. 
}
\ep
\end{remark}

 \vspace{3mm}

We next introduce the approximation error at time $n$ associated with the Hybrid Now algorithm by
\beqs
 \eps^{\mathrm{approx}}_{\mathrm{HN},n} &:=& \inf_{A \in \mathcal{A}_M}  \E_M\Big[ f\big( X_n, A(X_n)\big) +  \hat Y_{n+1}^{A} \Big] 
\; - \;  \inf_{A \in \A^{\Xc}} \E_M \Big[ f\big( X_n, A(X_n)\big) + \hat Y_{n+1}^{A} \Big],
\enqs
where  $\hat Y_{n+1}^{A}:=  \hat{V}_{n+1}^M \left( F\left( X_n, A(X_n), \eps_{n+1} \right) \right)$. 
We have the following bound on this approximation error:

\vspace{1mm}

\begin{lemma}
	\label{bound:approxErrorHybrid}
	For $n$ $=$ $0,...,N-1$, the following holds: 
	\begin{equation}
	\begin{split}
	\eps^{\mathrm{approx}}_{\mathrm{HN},n} &\leq \left( [f]_L + \lVert V_{n+1} \rVert_{\infty} [r]_L  \right)  \inf_{A \in \A^{\Xc}} \E_M \left[ \left\vert A(X_n) - a^{\mathrm{opt}}_n(X_n) \right\vert \right] \\
	& \hspace{0.5cm} + \;  2 \lVert r \rVert_{\infty} \E_M \left[ \left\vert V_{n+1}(X_{n+1}) - \hat{V}_{n+1}^M(X_{n+1}) \right\vert \right].
	\end{split}
	\end{equation}
\end{lemma}
\noindent \textbf{Proof.}  See Appendix \ref{sec:proofofLemmaApproxErrorHybrid}.
\ep

\vspace{3mm}

\noindent {\bf Proof of Theorem \ref{theoKohler}}  
%The proof is detailed in  Section \ref{sectionKoh} in Appendix.\\

\noindent  Observe that not only the optimal strategy but also the value function is estimated at each time step $n$ $=$ $0,...,N-1$ using neural networks in the hybrid algorithm. It spurs us to introduce the following auxiliary process $(\bar{V}_n^M)_{n=0}^{N}$ defined by backward induction as:
\begin{equation*}
\begin{cases}
\bar{V}_N^M(x)&=g(x), \quad \text{ for } x \in \Xc, \\
\bar{V}_n^M(x) &=   f(x,\hat a_n^M(x))+ \E \Big[  \hat{V}_{n+1}^M(F(x,\hat a_n^M(x),\eps_{n+1}))\Big], 
\quad \text{ for } x \in \Xc,
\end{cases}
\end{equation*} 
and we notice that for $n$ $=$ $0,...,N-1$,  $\bar{V}_n^M$  is  the quantity estimated by $\hat{V}_n^M$. 

\vspace{1mm}

\noindent {\it Step 1.} We state the following estimates: for $n$ $=$ $0,...,N-1$, 
	\begin{equation}
	\label{ineq:1_LemmaHybrid1}
		0 \leq \E_M \left[ \bar{V}_n^M(X_n) - \inf_{a\in A} \left\{ f(X_n,a) +  \E_{M,n,X_n}^a \left[ \hat V_{n+1}^M (X_{n+1}) \right]  \right\} \right] \leq 2 \eps^{\mathrm{esti}}_{\mathrm{HN},n} + \eps^{\mathrm{approx}}_{\mathrm{HN},n},
	\end{equation}
	and,
	\begin{align}
	\label{ineq:2_LemmaHybrid2}
	\E_M \left[ \left\vert \bar{V}_n^M(X_n) - \inf_{a\in A} \left\{ f(X_n,a) +  \E_{M,n,X_n}^a \left[ \hat V_{n+1}^M (X_{n+1}) \right]  \right\} \right \vert ^2 \right] \\
\leq 2\left((N-n) \lVert f \rVert_\infty +  \lVert g \rVert_\infty \right) \left( 2 \eps^{\mathrm{esti}}_{\mathrm{HN},n}+\eps^{\mathrm{approx}}_{\mathrm{HN},n} \right),\nonumber
\end{align}	
where $\E_{M,n,X_n}$ stands for the expectation conditioned by the training set and $X_n$.\\

\noindent Let us first show the estimate \eqref{ineq:1_LemmaHybrid1}. Note that  inequality 
\[ 
\bar{V}_n^M(X_n) - \inf_{a\in A} \left\{ f(X_n,a) +  \E_{M,n,X_n}^a \left[ \hat V_{n+1}^M (X_{n+1}) \right]  \right\}  \geq 0 
\] 
holds because $\hat{a}_n^M$ cannot do better than the optimal strategy. Take its expectation to get the first inequality in \eqref{ineq:1_LemmaHybrid1}. Moreover, we write 
\begin{align*}
\E_M \left[ \bar{V}_n^M(X_n) \right] &\leq	\E_M \left[ f\left( X_n, \hat a_n^M (X_n) \right) + \hat V_{n+1}^M \left( X_{n+1}^{\hat{a}_n^M} \right) \right] \\
	&\leq \inf_{A \in \mathcal{A}_M} \E_M \left[ f\left( X_n, A (X_n) \right) + \hat V_{n+1}^M \left( X_{n+1}^{A} \right) \right] + 2 \eps^{\mathrm{esti}}_{\mathrm{HN},n},
\end{align*}
which holds by the same arguments as those used to prove \reff{interestistep1}. We deduce that
\begin{align*}
\E_M \left[ \bar{V}_n^M(X_n) \right] 
%&\leq \E_M \left[ f\left( X_n, \hat A_n^M (X_n) \right) + \hat V_{n+1}^M \left( X_{n+1}^{\hat{A}_n^M} \right) \right]\\
&\leq \inf_{A \in \A^{\Xc}} \E_M \left[ f\left( X_n, A (X_n) \right) + \hat V_{n+1}^M \left( X_{n+1}^{A} \right) \right] + \eps^{\mathrm{approx}}_{\mathrm{HN},n}+2 \eps^{\mathrm{esti}}_{\mathrm{HN},n} \\
&\hspace{-.2cm}\leq \E_M \left[ \inf_{a \in A} \left\{   f\left( X_n, a\right) + \E_{M}^a \left[  \hat V_{n+1}^M \left( X_{n+1} \right) \big\vert X_n \right]    \right\}\right] + \eps^{\mathrm{approx}}_{\mathrm{HN},n}+2 \eps^{\mathrm{esti}}_{\mathrm{HN},n}. 
\end{align*}
This completes the  proof of the second inequality stated in \reff{ineq:1_LemmaHybrid1}.\\ On the other hand, noting:    
$$\left\vert \bar{V}_n^M(X_n) - \Inf_{a\in A} \left\{ f(X_n,a) +  \E_{M}^a \left[ \hat V_{n+1}^M (X_{n+1}) \big \vert X_n \right]  \right\} \right \vert  \leq 2 \left( (N-n) \lVert f \rVert_\infty +  \lVert g \rVert_\infty  \right),$$ and  using \reff{ineq:1_LemmaHybrid1}, we obtain the inequality   
\reff{ineq:2_LemmaHybrid2}. 

\vspace{2mm}

\noindent {\it Step 2.} We state the following estimation: for all $n \in \{0,...,N\}$
	\begin{align}
	\label{eq1Hybridinduction}
 & \left\lVert \hat{V}_{n}^{M}(X_{n})-\bar{V}_{n}^{M}(X_{n})\right\rVert _{M,1}\\
 & = \mathcal{O}_{\mathbb{P}}\Bigg(\gamma_{M}^{2}\sqrt{K_{M}\frac{log(M)}{M}}+\inf_{\Phi\in\mathcal{V}_{M}}\sqrt{\left\lVert \Phi(X_{n})-{V}_{n}^{M}(X_{n})\right\rVert _{M,1}}\nonumber\\
 & +\inf_{A\in\A^{\Xc}}\sqrt{\left\lVert A(X_{n})-a_{n}^{\mathrm{opt}}(X_{n})\right\rVert _{M,1}}+\sqrt{\left\lVert V_{n+1}(X_{n+1})-\hat{V}_{n+1}^{M}(X_{n+1})\right\rVert _{M,1}}\Bigg),\nonumber
\end{align}
	
	where ${\lVert . \rVert_{M,p} = \Big(\E_M \left[ \vert . \vert^p \right]}  \Big)^{1\over p}$, i.e. $\lVert . \rVert_{M,p}$ stands for the $\mathbb{L}^p$ norm conditioned by the training set,
	for $p \in \{1,2\}$. The proof relies on  Lemma \ref{lem:H1} and Lemma \ref{lem:H2} (see Appendix \ref{sectionKoh}).
 	 
Let us first show the following relation:
\begin{equation}
\label{eq1Hybridinduction_11}
\E_M\Big[ \big\vert  \hat{V}_n^M(X_n)-\bar{V}_n^M(X_n) \big\vert^2\Big] = \mathcal{O}_\mathbb{P} \bigg( \gamma_M^4K_M\frac{\log(M)}{M}+ \inf_{\Phi \in \mathcal{V}_M} \mathbb{E}\Big[ \vert \Phi(X_n)-\bar{V}_n^M(X_n) \vert^2 \Big] \bigg).
\end{equation}
For this, take $\delta_M=\gamma_M^4K_M\frac{\log(M)}{M}$, let $\delta > \delta_M$, and denote 
$$
\Omega_g:=\bigg\{ f-g: f \in \mathcal{V}_M, \frac{1}{M} \sum_{m=1}^{M} \big\vert f(x_m)- g(x_m) \big\vert^2 \leq \frac{\delta }{\gamma_M^2} \bigg\}.
$$ Apply Lemma \ref{lem:H2} to obtain:
\begin{align}
\label{ineq:lem11}
\int_{c_2\delta/\gamma_M^2}^{\sqrt{\delta}}  \log\left(  \mathcal{N}_2\bigg( \frac{u}{4\gamma_M}, \Omega_g, x^M_1 \bigg) \right)^{1/2} \diff u & \nonumber\\
&\hspace{-4.5cm} \leq \int_{c_2\delta/\gamma_M^2}^{\sqrt{\delta}}  \log\left(  \mathcal{N}_2\bigg( \frac{u}{4\gamma_M}, \mathcal{V}_M, x^M_1 \bigg) \right)^{1/2} \diff u \nonumber  \\
& \hspace{-4.5cm} \leq \int_{c_2\delta/\gamma_M^2}^{\sqrt{\delta}}   \big((4d+9)K_M+1\big)^{1/2} \left[  \log\left(\frac{48 e \gamma_M^2 \big( K_M +1\big)}{u}  \right) \right]^{1/2} \diff u  \nonumber \\
& \hspace{-4.5cm} \leq \int_{c_2\delta/\gamma_M^2}^{\sqrt{\delta}}   \big((4d+9)K_M+1\big)^{1/2} \left[  \log\left(48 e \frac{\gamma_M^4}{\delta} \big( K_M +1\big) \right) \right]^{1/2} \diff u \nonumber \\
& \hspace{-4.5cm} \leq \sqrt{\delta} \big((4d+9)K_M+1\big)^{1/2} \left[  \log\left(48 e \gamma_M^4M\big( K_M +1\big) \right) \right]^{1/2} \nonumber \\
& \hspace{-4.5cm} \leq c_5 \sqrt{\delta} \sqrt{K_M} \sqrt{ \log(M)},
\end{align}
where $\mathcal{N}_2(\eps, \mathcal{V},x_1^M)$ stands for the $\eps$-covering number of $\mathcal{V}$ on $x_1^M$, which is introduced in section \ref{sectionKoh}, and where the last line holds since we assumed $\frac{M \delta_M}{\gamma_M^2} \xrightarrow[M \to 0]{} 0$.  Since $\delta$ $>$ $\delta_M:=\gamma_M^4K_M\frac{\log(M)}{M}$, we then have 
$\sqrt{\delta} \sqrt{K_M} \sqrt{ \log(M)} \leq \frac{\sqrt{M} \delta}{\gamma_M^2}$, which implies that \eqref{ineq:condtionCoveringNumber} holds by \eqref{ineq:lem11}. 
Therefore, by application of Lemma \ref{lem:H1}, the following holds:
\begin{equation*}
\label{eq1Hybridinduction_2}
\E_M\Big[ \big\vert  \tilde{V}_n^M(X_n)-\bar{V}_n^M(X_n) \big\vert^2\Big] = \mathcal{O}_\mathbb{P} \bigg( \gamma_M^4K_M\frac{\log(M)}{M}+ \inf_{\Phi \in \mathcal{V}_M} \mathbb{E}\Big[ \vert \Phi(X_n)-\bar{V}_n^M(X_n) \vert^2 \Big] \bigg).
\end{equation*}
It remains to note that $ \E_M\Big[ \big\vert  \hat{V}_n^M(X_n)-\bar{V}_n^M(X_n) \big\vert^2\Big] \leq \E_M\Big[ \big\vert  \tilde{V}_n^M(X_n)-\bar{V}_n^M(X_n) \big\vert^2\Big] $ always holds, and this completes the proof of \reff{eq1Hybridinduction_11}.

Next, let us show 
\beq 
 & & \hspace{-1cm} \inf_{\Phi \in \mathcal{V}_M} \left \lVert \Phi(X_n)- \bar V_n(X_n) \right \rVert_{M,2} \nonumber  \\
 &\hspace{-.5cm}=&  \Oc \Bigg( \gamma_M^2 \sqrt{\frac{K_M \log(M)}{M}} + \sup_{n \leq k \leq N} \inf_{\Phi \in \mathcal{V}_M} \left \lVert \Phi(X_n) - V_n(X_n)  \right \rVert_{M,2} \nonumber   \\
& &   + \;  \inf_{A \in \mathcal{A}_M} \E_M \left[ \left \vert A(X_n)- a^{\mathrm{opt}}_n(X_n) \right\vert \right] +\left  \lVert V_{n+1}(X_{n+1}) - \hat{V}_{n+1}^M(X_{n+1})\right \rVert_{M,2} \Bigg).  
\label{estimstep2}
\enq
For this,  take some arbitrary $\Phi \in \mathcal{V}_M$ and split
% \begin{align*}
%\left \lVert \Phi(X_n)- \bar V_n^M(X_n) \right \rVert_{M,2} \leq \left \lVert \Phi(X_n)- V_n(X_n) \right \rVert_{M,2} + \left \lVert V_n(X_n)- \bar V_n^M(X_n) \right \rVert_{M,2},
%\end{align*}
%so that:
\begin{align}
\label{ineq:splitHybrid1}
 \inf_{\Phi \in \mathcal{V}_M}	\left \lVert \Phi(X_n)- \bar V_n^M(X_n) \right \rVert_{M,2} \;& \leq \; 
  \inf_{\Phi \in \mathcal{V}_M}\left \lVert \Phi(X_n)- V_n(X_n) \right \rVert_{M,2} \\
   &\hspace{.4cm}\; + \left \lVert V_n(X_n)- \bar V_n^M(X_n) \right \rVert_{M,2}.
\end{align}
To bound the last term in the r.h.s. of \eqref{ineq:splitHybrid1}, we write
 \begin{align*}
\left \lVert V_n(X_n)- \bar V_n^M(X_n) \right \rVert_{M,2} & \\
&\hspace{-2cm}\leq \left \lVert V_n(X_n)- \inf_{a \in A} \left\{   f\left( X_n, a\right) + \E_{M}^a \left[  \hat V_{n+1}^M \left( X_{n+1} \right) \big\vert X_n \right]    \right\}\right \rVert_{M,2} \\
& \hspace{-1.7cm}+ \left \lVert\inf_{a \in A} \left\{   f\left( X_n, a\right) + \E_{M}^a \left[  \hat V_{n+1}^M \left( X_{n+1} \right) \big\vert X_n \right]    \right\} - \bar V_n^M(X_n) \right \rVert_{M,2}. 
\end{align*}
Use the dynamic programming principle, assumption \ref{Hd} and \eqref{ineq:2_LemmaHybrid2} to get:
 \begin{align*}
\left \lVert V_n(X_n)- \bar V_n^M(X_n) \right \rVert_{M,2} &\leq \lVert r \rVert_{\infty} \left \lVert V_{n+1}(X_{n+1})-  \hat V_{n+1}^M \left( X_{n+1} \right)\right \rVert_{M,2} \\
& \hspace{.5cm}+  \sqrt{2\left((N-n) \lVert f \rVert_\infty +  \lVert g \rVert_\infty \right) \left( 2 \eps^{\mathrm{esti}}_{\mathrm{HN},n} + \eps^{\mathrm{approx}}_{\mathrm{HN},n} \right) }.
\end{align*}
%\E_M \left[ \inf_{a \in A} \left\{   f\left( X_n, a\right) + \E_{M}^a \left[  \hat V_{n+1}^M \left( X_{n+1} \right) \big\vert X_n \right]    \right\}\right]
We then notice that  
\begin{align}
	\left\vert V_{n+1}(X_{n+1})-\hat V_{n+1}^M(X_{n+1}) \right\vert^2 & \nonumber\\
	& \hspace{-2cm} \leq 2 \lVert r \rVert_{\infty} \left((N-n) \lVert f \rVert_\infty +  \lVert g \rVert_\infty \right) \left\vert V_{n+1}(X_{n+1})-\hat V_{n+1}^M(X_{n+1}) \right\vert
\end{align}
 holds a.s., so that
 \begin{align*}
\left \lVert V_n(X_n)- \bar V_n^M(X_n) \right \rVert_{M,2} &\\
& \hspace{-2cm}\leq \sqrt{ 2 \lVert r \rVert_{\infty} \left((N-n) \lVert f \rVert_\infty +  \lVert g \rVert_\infty \right) \left \lVert V_{n+1}(X_{n+1})-  \hat V_{n+1}^M \left( X_{n+1} \right)\right \rVert_{M,1} } \\
& \hspace{-1.5cm}+  \sqrt{2\left((N-n) \lVert f \rVert_\infty +  \lVert g \rVert_\infty \right) \left( 2 \eps^{\mathrm{esti}}_{\mathrm{HN},n} + \eps^{\mathrm{approx}}_{\mathrm{HN},n} \right) },
\end{align*}
and  use Lemma \ref{bound:approxErrorHybrid} to bound $\eps^{\mathrm{approx}}_{\mathrm{HN},n}$.
By plugging into \reff{ineq:splitHybrid1}, and using the estimations in Lemmas  \ref{speedLearning_Hybrid} and \ref{bound:approxErrorHybrid}, we obtain the  estimate \reff{estimstep2}. 
Together with \reff{eq1Hybridinduction_11}, this proves the required estimate\reff{eq1Hybridinduction}. 
 By induction, we get as $M \to \infty$,
	\beqs
	\E_M\Big[\vert  \hat{V}_n^M(X_n)-V_n(X_n) \vert  \Big]&=&  \mathcal{O}_\mathbb{P}\Bigg(  \left( \gamma_M^4 \frac{K_M\log(M)}{M} \right)^{\frac{1}{2(N-n)}} + \left( \gamma_M^4 \frac{\rho_M^2 \eta_M^2 \log(M)}{M} \right)^{\frac{1}{4(N-n)}} \\
	& & \;\;\;\;\; + \;  \sup_{n \leq k \leq N} \inf_{\Phi \in \mathcal{V}_M} \left(\E_M\Big[ \vert \Phi(X_k)-V_k(X_k) \vert^2 \Big] \right)^{\frac{1}{2(N-n)}}  \nonumber \\
	 & & \;\;\;\;\; + \;   \sup_{n \leq k \leq N}\inf_{A \in \mathcal{A}_M} \left(  \mathbb{E}\Big[ \vert A(X_k)-a^{\mathrm{opt}}_k(X_k) \vert \Big] \right)^{\frac{1}{2(N-n)}}\Bigg),
	\enqs
	which completes the proof of Theorem \ref{theoKohler}. 
\ep

\section{Conclusion} \label{secconclusion}
This paper develops new machine learning algorithms for high-dimensional  Markov decision processes. We propose and compare two algorithms based on dynamic programming, performance/hybrid iteration and neural networks for approximating the control and value function. The main theoretical contribution is to provide a detailed convergence analysis for each of these algorithms: by using least squares neural network regression,  we state error estimates in terms of the universal approximation error of  neural networks, and of the statistical risk when estimating the network function.  Numerical tests on various applications are presented in a companion paper \cite{bacetal18}.

\appendix

\section{Localization} \label{sec::ProblemLocalization}

In this section, we show how to relax the boundedness condition on  the state space by a localization argument.

Let $R>0$. Consider the localized state space $\bar B_d^\Xc(0,R):=\Xc \cap \{ \lVert x \rVert_d \leq R \}$, where $\lVert . \rVert_d$  is the Euclidean norm of $\R^d$. Let $\left( \bar X_n \right)_{0 \leq n \leq N}$ be the Markov chain defined by its transition probabilities as
\[
\P \left( \bar{X}_{n+1} \in O \Big\vert \bar X_n=x,a \right) = \int_O r(x,a;y) d\pi_R\circ \mu(y), \quad \text{ for } n=0,\ldots, N-1,
\]
for all Borelian $O$ in $\bar B_d^\Xc(0,R)$, where $\pi_R$ is the Euclidean projection of $\R^d$ on $\bar B_d^\Xc(0,R)$, and $\pi_R\circ \mu$ is the pushforward measure of $\mu$. 
Notice that the transition probability of $\bar X$ admits the same density $r$, for which \ref{Hd} holds, w.r.t. $\pi_R \circ \mu$.

\vspace{1mm}

\noindent Define $(\bar{V}^R_n)_n$ as the value function associated with the following stochastic control problem for $\left( \bar X_n \right)_{n=0}^{N}$:
\begin{equation}
\label{def:VBar}
\left\{  
\begin{array}{ccl}
\bar V_N^R(x) &= &g(x), \qquad \\
\bar V_n^R(x) & = &\Inf_{ \alpha \in \mathcal{C}} \; \E \left[ \sum_{k=n}^{N-1} f\left( \bar X_k,\alpha_k \right) + g \left(\bar X_N \right)  \right], \text{ for } n=0,\ldots,N-1,\\
\end{array}  
 \right.
\end{equation}
 for  $x \in \bar{B}_d^\Xc(0,R)$.  By the dynamic programming principle, $(\bar{V}^R_n)_n$ is solution of the following Bellman backward equation: 
\begin{equation}
\begin{cases}
\bar{V}^R_N(x)=g(x)\\
\bar{V}^R_n(x)=\Inf_{a \in \A} \; \bigg\{f(x,a)+ \mathbb{E}_n^a\Big[\bar{V}^R_{n+1}\big(\pi_R\big(F(x,a,\eps_{n+1}) \big) \big) \Big] \bigg\}, &\quad \forall x \in B_d^\Xc(0,R). \nonumber
\end{cases}
\end{equation}

\noindent We shall assume two conditions on the measure $\mu$.\\
\noindent \setword{\textbf{(Hloc)}}{Hloc} $\mu$ is such that:
	\begin{equation*}
	\mathbb{E} \big[ \vert \pi_R(X)- X \vert \big] \xrightarrow[R\to \infty]{} 0\quad \text{ and } \quad \P \left(\vert X \vert >R \right)  \xrightarrow[R\to \infty]{} 0, \quad \text{ where $X \sim \mu$. }
	\end{equation*}
Using the dominated convergence theorem, it is straightforward to see that \ref{Hloc} holds if  $\mu$ admits a moment of order 1.

\begin{proposition}
	\label{prop:localize}
	Let $X_n \sim \mu$. The following holds:
	\begin{align}
	\E\Big[ \Big\vert \bar{V}^R_n \big(\pi_R(X_n)\big) - V_n\big(X_n\big) \Big\vert \Big] & \nonumber\\
	&\hspace{-2cm}\leq   \lVert V \rVert_{\infty} \Big( [r]_L \E \left[ \vert \pi_R(X_n) - X_n\vert \right]  +2 \P \left( \vert X_n \vert > R\right) \Big) \frac{1- \lVert r \rVert_{\infty}^{N-n}}{1- \lVert r \rVert_{\infty}} \nonumber\\
	& \hspace{-1.5cm} +[g]_L  \lVert r \rVert_{\infty}^{N-n} \E \left[ \vert \pi_R(X_n) - X_n\vert \right], \nonumber
	\end{align}
	where we denote $\lVert V \rVert_{\infty}=\Sup_{0 \leq k \leq N} \; \lVert V_{k} \rVert_{\infty}$, and use the convention $\frac{1-x^p}{1-x}=p$ for $x=0$ and $p>1$. 
	Consequently, for all $n=0,...,N$, we get under \ref{Hloc}:
	\begin{equation}
		\E\Big[ \Big\vert \bar{V}^R_n \big(\pi_R(X_n)\big) - V_n\big(X_n\big) \Big\vert \Big] \xrightarrow[R \to \infty]{}0, \quad \text{ where } X_n \sim \mu. \nonumber
	\end{equation}
\end{proposition}

\noindent \textbf{Comment: } Proposition \ref{prop:localize} states that if $\Xc$ is not bounded, the control problem \eqref{def:VBar} associated with a bounded controlled process $\bar X$  can be as close as desired, in $\mathbb{L}^1(\mu)$ sense, to the original control problem by taking $R$ large enough. Moreover, as stated before, the transition probability of $\bar X$ admits the same density $r$ as $X$ w.r.t. the pushforward measure $\pi_R \circ \mu$. 

\vspace{2mm}

\noindent \textbf{Proof of Proposition \ref{prop:localize}.}  
Take $X_n \sim \mu$ and write:
\begin{align}
\label{ineq:localization}
\mathbb{E}\Big[ \big\vert   \bar{V}^R_n(\pi_R(X_n))-V_n(X_n) \big\vert \Big]  &\leq \E \left[   \big\vert   \bar{V}^R_n(X_n)-V_n(X_n) \big\vert \mathds{1}_{\vert X_n\vert \leq R}\right] \nonumber  \\
& \hspace{.5cm} +  \E \left[   \big\vert   \bar{V}^R_n(\pi_R(X_n))-V_n(X_n) \big\vert \mathds{1}_{\vert X_n\vert \geq R}\right].
\end{align}

\noindent Let us first bound the first term in the r.h.s. of \eqref{ineq:localization}, by showing that, for $n=0,\ldots, N-1$:
\begin{align}
	 \E \left[   \big\vert   \bar{V}^R_n(X_n)-V_n(X_n) \big\vert \mathds{1}_{\vert X_n\vert \leq R}\right] &\leq \lVert r \rVert_{\infty} \mathbb{E}\Big[ \big\vert   \bar{V}^R_{n+1}(\pi_R(X_{n+1}))-V_{n+1}(X_{n+1}) \big\vert \Big]  \nonumber \\
	 & \hspace{-1cm} + [r]_L \lVert V_{n+1} \rVert_{\infty} \E \left[ \vert \pi_R(X_{n+1}) - X_{n+1} \vert \right], \text{ with } X_{n+1} \sim \mu. \label{step1} 
\end{align}
Take $x \in \bar B_d(0,R)$ and notice that
\begin{align*}
	\big\vert   \bar{V}^R_n(x)-V_n(x) \big\vert &\leq \inf_{a \in A}  \bigg\{ \int_A \left\vert \bar V_{n+1}^R\left( \pi_R(y)  \right) - V_{n+1}(y) \right\vert r\left(x,a;\pi_R(y) \right) d \mu(y) \\
	&\hspace{1.5cm} + \int \left\vert V_{n+1}(y) \right\vert \left\vert r(x,a; \pi_R(y) ) - r(x,a;y) \right\vert d \mu(y) \bigg\} \\
	& \leq \lVert r \rVert_{\infty} \E\left[  \left\vert \bar V_{n+1}^R\left( \pi(X_{n+1}) \right) - V_{n+1}\left(X_{n+1}\right) \right\vert \right]\\
	& \hspace{.5cm} + [r]_L \lVert V_{n+1} \rVert_{\infty} \E \left[  \left\vert \pi_R(X_{n+1})-  X_{n+1} \right\vert \right], \quad \text{ where } X_{n+1} \sim \mu.
\end{align*}
It remains to inject this bound in the expectation to obtain \reff{step1}.

To bound the second term in the r.h.s. of \eqref{ineq:localization}, notice that  
\[
\big\vert   \bar{V}^R_n(\pi_R(X_n))-V_n(X_n) \big\vert \leq 2 \lVert V_n \rVert_{\infty}
\]  holds a.s., which implies:
\begin{equation} \label{step2}
	\E \left[   \big\vert   \bar{V}^R_n(\pi_R(X_n))-V_n(X_n) \big\vert \mathds{1}_{\vert X_n\vert \geq R}\right] \leq 2 \lVert V_n \rVert_{\infty} \P \left( \vert X_n \vert > R\right).
\end{equation}

\vspace{1mm}

\noindent Plugging \reff{step1} and \reff{step2} into \eqref{ineq:localization} yields:
\begin{align}
	\mathbb{E}\Big[ \big\vert   \bar{V}^R_n(\pi_R(X_n))-V_n(X_n) \big\vert \Big]  &\leq \lVert r \rVert_{\infty} \mathbb{E}\Big[ \big\vert   \bar{V}^R_{n+1}(\pi_R(X_{n+1}))-V_{n+1}(X_{n+1}) \big\vert \Big] \nonumber  \\
	& \hspace{-1.5cm} + [r]_L \lVert V_{n+1} \rVert_{\infty} \E \left[ \vert \pi_R(X_{n+1}) - X_{n+1} \vert \right]  +2 \lVert V_n \rVert_{\infty} \P \left( \vert X_n \vert > R\right), \nonumber
\end{align}
with $X_n$ and $X_{n+1}$ i.i.d. following the law $\mu$.  The result stated in proposition \ref{prop:localize} then follows by induction.
\ep

\section{Forward evaluation of the optimal controls in $\mathcal{A}_M$}

\label{sec::forwardEvaluation}

We evaluate in this section the real performance of the best controls in $\mathcal{A}_M$. 
Let $(a^{\mathcal{A}_M}_n)_{n=0}^{N-1}$  be the sequence of optimal controls in the class of neural networks $\mathcal{A}_M$, and denote by $(J^{\mathcal{A}_M}_n)_{0\leq n \leq N}$ the cost functional sequence associated with $(a^{\mathcal{A}_M}_n)_{n=0}^{N-1}$ and characterized as solution of the Bellman equation:
	\[\begin{cases}
	J^{\mathcal{A}_M}_N(x)= g(x) \\
	J^{\mathcal{A}_M}_n(x) = \Inf_{A \in \mathcal{A}_M} \left\{ f(x, A(x)) + \mathbb{E}_{n,X_n}^A [J^{\mathcal{A}_M}_{n+1}(X_{n+1})]  \right\},
	\end{cases}\]
	where $\mathbb{E}_{n,X_n}^A[\cdot]$ stands for the expectation conditioned by $X_n$ and when the control $A$ is applied at time $n$.
	
	In this section, we are interested in comparing  $J^{\mathcal{A}_M}_n$ to $V_n$. 
	Note that $V_n(x) \leq J^{\mathcal{A}_M}_n(x)$ holds for all $x \in \Xc$, since $\mathcal{A}_M$ is included in the set of the Borelian functions of $\Xc$. We can actually show the following:

	\begin{proposition}
		\label{propForward}
		Assume that there exists a sequence of optimal feedback controls $(a^{\mathrm{opt}}_n)_{0 \leq n\leq N-1}$ for the control problem with value function $V_n$, $n=0,\ldots,N$.
		Then the following holds, as $M \to \infty$:
		\begin{equation}
		\label{eq::propForward}
			\E \left[ J^{\mathcal{A}_M}_n(X_n)- V_n(X_n) \right]  = \mathcal{O} \left( \sup_{n \leq k \leq N-1}\inf_{A \in \mathcal{A}_M} \E \left[  \lvert A(X_k) - a_k^{\mathrm{opt}}(X_k) \rvert \right] \right).
		\end{equation}
	\end{proposition}

%	Proposition \ref{ConvergenceError} shows that the quantity in the r.h.s. of \eqref{eq::propForward} cancels when $M \to +\infty$, under two different %regularity assumptions on the optimal control.
\begin{remark}
{\rm
	Notice  that there is no estimation error term in \eqref{eq::propForward}, since the optimal strategies in $\mathcal{A}_M$ are defined as those minimizing the real cost functionals in 
	$\mathcal{A}_M$, and not the empirical ones.
	}
\ep	
\end{remark}

\vspace{2mm}

\noindent {\bf Proof of Proposition \ref{propForward}.}  
	Let $n \in \{0,...,N-1\}$, and $X_n \sim \mu$. Take $A \in \mathcal{A}_M$, and denote $J^{A}_n(X_n) = f(x, A(x)) + \mathbb{E}_{n,X_n}^A [J^{\mathcal{A}_M}_{n+1}(X_{n+1})] $. 
	 Clearly, we have $J^{\mathcal{A}_M}_{n}= \Min_{A \in \mathcal{A}_M} J^{A}_{n}$. Moreover:
	\begin{align}
	\E \Big[ J^{A}_n(X_n)- V_n(X_n) \Big] & \leq  \E \Big[  \lvert f(X_n,A(X_n)) - f(X_n, a^{\mathrm{opt}}_n(X_n)) \rvert \Big]  \nonumber \\
	& \hspace{-2cm}+ \E \Big[  \lvert J^{\mathcal{A}_M}_{n+1}(F(X_n,A(X_n),\eps_{n+1})) -   V_{n+1}(F(X_n,a^{\mathrm{opt}}_n(X_n),\eps_{n+1})) \rvert \Big]  \nonumber\\
	& \leq [f]_L  \E \Big[  \lvert a^{\mathrm{opt}}_n(X_n) - A(X_n) \rvert \Big]  \nonumber \\
	& \hspace{-2cm} +   \E \Big[ \lvert V_{n+1}(F(X_n,A(X_n),\eps_{n+1})) -   V_{n+1}(F(X_n,a^{\mathrm{opt}}_n(X_n),\eps_{n+1})) \rvert \Big] \nonumber \\
	&\hspace{-2cm} +   \E\Big[ \lvert J^{\mathcal{A}_M}_{n+1}(F(X_n,A(X_n),\eps_{n+1})) -   V_{n+1}(F(X_n,A(X_n),\eps_{n+1})) \rvert \Big].	\label{ineq::approxerror}
	\end{align}
Applying assumption \ref{Hd} to bound the last term in the r.h.s. of \eqref{ineq::approxerror} yields
	\begin{align}
	 \E \Big[ J^{A}_n(X_n)- V_n(X_n) \Big] & \leq  \big( [f]_L + \lVert V_{n+1} \rVert_{\infty} [r]_L\big) \E \Big[  \lvert a^{\mathrm{opt}}_n(X_n) - A(X_n) \rvert \Big] \nonumber\\
	 &\hspace{.5cm} + \lVert r \rVert_\infty  \E \Big[ \lvert J^{\mathcal{A}_M}_{n+1} (X_{n+1}) -V_{n+1}(X_{n+1}) \rvert \Big],  \nonumber
	\end{align}
	which holds for all $A \in \mathcal{A}_M$, so that:
	\begin{align}
\E \Big[  J^{\mathcal{A}_M}_n(X_n)- V_n(X_n) \Big] & \leq  \big( [f]_L + \lVert V_{n+1} \rVert_{\infty} [r]_L\big)  \Inf_{A \in \mathcal{A}_M} \E \Big[  \lvert a^{\mathrm{opt}}_n(X_n) - A(X_n) \rvert \Big] \nonumber\\
&\hspace{.5cm} + \lVert r \rVert_\infty  \E \Big[ \lvert J^{\mathcal{A}_M}_{n+1} (X_{n+1}) -V_{n+1}(X_{n+1}) \rvert \Big].  \nonumber
\end{align}	
	\eqref{eq::propForward} then follows directly by induction. 
\ep

\section{Proof of Lemma \ref{speedLearning}}
\label{sec::estimationError}

The proof is  divided into four steps.

\noindent \textit{Step 1: Symmetrization by a ghost sample.} We take $\eps >0$ and show that for\\ ${M > 2\frac{\big( (N-n) \lVert f \rVert_\infty + \lVert g \rVert_\infty \big)^2}{\eps^2}}$, the following holds:
{\small 
	\begin{align}
 & \mathbb{P}\Bigg[\sup_{A\in\mathcal{A}_{M}}\Bigg\vert\frac{1}{M}\sum_{m=1}^{M}\Big[f(X_{n}^{(m)},A(X_{n}^{(m)}))+\hat{Y}_{n+1}^{(m),A}\Big]-\mathbb{E}\Big[J_{n}^{A,(\hat{a}_{k}^{M})_{k=n+1}^{N-1}}(X_{n})\Big]\Bigg\vert>\eps\Bigg]\nonumber\\
 &\leq2\mathbb{P}\Bigg[\!\!\sup_{A\in\mathcal{A}_{M}}\!\Bigg\vert\frac{1}{M}\!\sum_{m=1}^{M}\!\Big[f(X_{n}^{(m)}\!\!,A(X_{n}^{(m)}))+\hat{Y}_{n+1}^{(m),A}\!\!-f(X_{n}^{'(m)}\!\!,A(X_{n}^{'(m)}))-\hat{Y}_{n+1}^{\prime\,(m),A}\Big]\Bigg\vert\!>\!\frac{\eps}{2}\Bigg]\!,
 \label{ineq::step1}
\end{align}}
where:
\begin{itemize}
	\item $(X^{'(m)}_k)_{1\leq m \leq M, n\leq k\leq N}$ is a copy of $(X^{(m)}_k)_{1\leq m \leq M, n\leq k\leq N}$ generated from an independent copy of the exogenous noises 
	$(\eps'^{(m)}_k)_{1\leq m \leq M, n\leq k\leq N}$, and independent copy of initial positions at time $n$, $(X^{'(m)}_n)_{m=1}^{M}$, following the same control $\hat{a}^M_k$ at time $k$$=n+1, \dots,  N-1$, and control $A$ at time $n$,
	\item We recall that $Y^{(m),A}_{n+1}$ has already been defined in \eqref{defY}, and we similarly define 
	\begin{equation*}
	Y^{\prime \, (m),A}_{n+1}:=  \sum_{k=n+1}^{N-1} f(X_k^{\prime \, (m),A},\hat a_k^M(X_k^{\prime \,  (m),A})) +  g(X_N^{\prime \, (m),A}).
	\end{equation*}
\end{itemize} 
Let $A^* \in \mathcal{A}_M$ be such that:
\begin{equation}
\Bigg\vert \frac{1}{M} \sum_{m=1}^M \Big[ f(X_n^{(m)},A^*(X_n^{(m)}))  +\hat Y_{n+1}^{(m),A^*} \Big] - \mathbb{E}\Big[ J_n^{A^*,(\hat{a}^M_k)_{k=n+1}^{N-1}}(X_n) \Big] \Bigg\vert > \eps \nonumber
\end{equation}
	if such a function exists, and an arbitrary function in $\mathcal{A}_M$ if such a function does not exist. Note that $\frac{1}{M} \sum_{m=1}^M \Big[ f(X_n^{(m)},A^*(X_n^{(m)}))  +\hat Y_{n+1}^{(m),A^*} \Big] - \mathbb{E}\Big[ J_n^{A^*,(\hat{a}^M_k)_{k=n+1}^{N-1}}(X_n) \Big] $ is a r.v., which implies that $A^*$  also depends on $\omega \in \Omega$.
Denote by $\P_M$ the probability conditioned by the training set of exogenous noises $(\eps^{(m)}_k)_{1 \leq m \leq M, n \leq k \leq N}$ and initial positions $(X^{(m)}_k)_{1\leq m \leq M, n\leq k\leq N}$, and recall  that $\E_M$ stands for the expectation conditioned by the latter. An application of Chebyshev's inequality yields
\begin{align}
\P_M \Bigg[  \left\vert \E_M\Big[ J_n^{A^*,(\hat{a}^M_k)_{k=n+1}^{N-1}}(X'_n) \Big] -  \frac{1}{M} \sum_{m=1}^M  \Big[ f(X'^{(m)}_n,A^*(X'^{(m)}_n))  +\hat Y_{n+1}^{\prime \, (m),A^*} \Big] \right\vert > \frac{\eps}{2} \Bigg]& \nonumber\\
& \hspace{-12cm}\leq \frac{\mathrm{Var}_M\Big[  J_n^{A^*,(\hat{a}^M_k)_{k=n+1}^{N-1}}(X'_n) \Big] }{M(\eps/2)^2}  \leq  \frac{\big( (N-n) \lVert f \rVert_\infty + \lVert g \rVert_\infty \big)^2}{M  \eps^2}, \nonumber
\end{align} 
where we have used $ 0 \leq \Big\vert J_n^{A^*,(\hat{a}^M_k)_{k=n+1}^{N-1}}(X'_n) \Big\vert  \leq (N-n) \lVert f \rVert_\infty + \lVert g \rVert_\infty$ which implies
{\small \begin{align*}
 & \hspace{-4mm}\mathrm{Var}_{M}\left[J_{n}^{A^{*},(\hat{a}_{k}^{M})_{k=n+1}^{N-1}}(X'_{n})\right]=\mathrm{Var}_{M}\bigg[J_{n}^{A^{*},(\hat{a}_{k}^{M})_{k=n+1}^{N-1}}(X'_{n})-\frac{(N-n)\lVert f\rVert_{\infty}+\lVert g\rVert_{\infty}}{2}\bigg]\\
 & \leq\E\left[\left(J_{n}^{A^{*},(\hat{a}_{k}^{M})_{k=n+1}^{N-1}}(X'_{n})-\frac{(N-n)\lVert f\rVert_{\infty}+\lVert g\rVert_{\infty}}{2}\right)^{2}\right]\leq\frac{\big((N-n)\lVert f\rVert_{\infty}+\lVert g\rVert_{\infty}\big)^{2}}{4}.
\end{align*}}
Thus, for $M > 2\frac{\big( (N-n) \lVert f \rVert_\infty + \lVert g \rVert_\infty \big)^2}{\eps^2}$, we have
\begin{align}
\label{ineq::step1_1}
\P_M \Bigg[  \left\vert \E_M\Big[ J_n^{A^*,(\hat{a}^M_k)_{k=n+1}^{N-1}}(X_n) \Big] -  \frac{1}{M} \sum_{m=1}^M  \Big[ f(X'^{(m)}_n,A^*(X'^{(m)}_n))  +\hat Y_{n+1}^{\prime \, (m),A^*} \Big]  \right\vert \leq \frac{\eps}{2} \Bigg] & \nonumber \\
& \hspace{-2cm} \geq \frac{1}{2}.
\end{align}
Hence:
{\small 
	\begin{align*}
\mathbb{P}\Bigg[\sup_{A \in \mathcal{A}_M} \Bigg\vert \frac{1}{M} \sum_{m=1}^M \Big[ f(X_n^{(m)},A(X_n^{(m)}))  +\hat Y_{n+1}^{(m),A} - f(X'^{(m)}_n,A(X'^m_n))  -\hat Y^{\prime \, (m),A}_{n+1} \Big]   \Bigg\vert  > \frac{\eps}{2}\Bigg] \\
& \hspace{-12cm} \geq  	\mathbb{P}\Bigg[ \Bigg\vert \frac{1}{M} \sum_{m=1}^M \Big[ f(X_n^{(m)},A^*(X_n^{(m)}))  +\hat Y_{n+1}^{(m),A^*}  - f(X'^{(m)}_n,A^*(X'^{(m)}_n))  -\hat Y^{\prime \, (m),A^*}_{n+1} \Big]   \Bigg\vert  > \frac{\eps}{2}\Bigg] \\
& \hspace{-12cm} \geq 	\mathbb{P}\Bigg[ \Bigg\vert \frac{1}{M} \sum_{m=1}^M \Big[ f(X_n^{(m)},A^*(X_n^{(m)}))  +\hat Y_{n+1}^{(m),A^*} \Big] - \E_M\Big[ J_n^{A^*,(\hat{a}^M_k)_{k=n+1}^{N-1}}(X_n) \Big] \Bigg\vert > \eps,\\
&  \hspace{-11.1cm} \Bigg\vert \frac{1}{M} \sum_{m=1}^M \Big[ f(X'^{(m)}_n,A^*(X'^{(m)}_n))  +\hat Y^{\prime \, (m),A^*}_{n+1} \Big]  - \E_M\Big[ J_n^{A^*,(\hat{a}^M_k)_{k=n+1}^{N-1}}(X_n) \Big]  \Bigg\vert  \leq  \frac{\eps}{2}\Bigg].
\end{align*}}
Observe that $\frac{1}{M} \sum_{m=1}^M \Big[ f(X_n^{(m)},A^*(X_n^{(m)}))  +\hat Y^{(m),A^*}_{n+1} \Big] - \E_M\Big[ J_n^{A^*,(\hat{a}^M_k)_{k=n+1}^{N-1}}(X_n) \Big] $ is measurable w.r.t. the $\sigma$-algebra generated by the training set, so that conditioning by the training set and injecting \eqref{ineq::step1_1} yields
{\small
\begin{align*}
2 \mathbb{P}\Bigg[\sup_{A \in \mathcal{A}_M} \Bigg\vert \frac{1}{M} \sum_{m=1}^M \Big[ f(X_n^{(m)},A(X_n^{(m)}))  +\hat Y^{(m),A}_{n+1} - f(X'^{(m)}_n,A(X'^{(m)}_n))  -\hat Y^{\prime \, (m),A}_{n+1} \Big]   \Bigg\vert  > \frac{\eps}{2}\Bigg] &\\
& \hspace{-12cm} \geq \mathbb{P}\Bigg[ \Bigg\vert \frac{1}{M} \sum_{m=1}^M \Big[ f(X_n^{(m)},A^*(X_n^{(m)}))  +\hat Y^{(m),A^*}_{n+1} \Big] - \E_M\Big[ J_n^{A^*,(\hat{a}^M_k)_{k=n+1}^{N-1}}(X_n) \Big] \Bigg\vert > \eps \Bigg]\\
& \hspace{-12cm} = \mathbb{P}\Bigg[\sup_{A \in \mathcal{A}_M} \Bigg\vert \frac{1}{M} \sum_{m=1}^M \Big[ f(X_n^{(m)},A(X_n^{(m)}))  +\hat Y^{(m),A}_{n+1} \Big] - \mathbb{E}\Big[ J_n^{A,(\hat a_k^M)_{k=n+1}^{N-1}}(X_n) \Big] \Bigg\vert > \eps \Bigg]
\end{align*}}
for $M > 2\frac{\big( (N-n) \lVert f \rVert_\infty + \lVert g \rVert_\infty \big)^2}{\eps^2}$, where we use the definition of $A^*$ to go from the second-to-last to the last line. The proof of \eqref{ineq::step1} is then completed. 

\vspace{2mm}

\noindent \textit{Step 2:}  We show that 
\begin{align}
\mathbb{E}\Bigg[\sup_{A \in \mathcal{A}_M} \Bigg\vert \frac{1}{M} \sum_{m=1}^M \Big[ f(X_n^{(m)},A(X_n^{(m)}))  +\hat Y^{(m),A}_{n+1} \Big] - \mathbb{E}\Big[ J_n^{A,(\hat a_k^M)_{k=n+1}^{N-1}}(X_n) \Big] \Bigg\vert \Bigg] & \nonumber \\
&\hspace{-12cm} \leq 4\mathbb{E}\Bigg[\sup_{A \in \mathcal{A}_M} \Bigg\vert \frac{1}{M} \sum_{m=1}^M \Big[ f(X_n^{(m)},A(X_n^{(m)}))  +\hat Y^{(m),A}_{n+1} - f(X'^{(m)}_n,A(X'^{(m)}_n))  
-\hat Y^{\prime \, (m),A}_{n+1} \Big]   \Bigg\vert \Bigg] \nonumber \\
& \hspace{-9.5cm} + \mathcal{O}\left(  \frac{1}{\sqrt{M}}  \right). \label{ineqstep2}
\end{align}
Indeed, let $M'=\sqrt{2} \frac{(N-n) \lVert f \rVert_\infty + \lVert g \rVert_\infty}{\sqrt{M}}$, and notice
{\small
	\begin{align}
 \label{ineq::interm333}
\mathbb{E}\Bigg[\sup_{A \in \mathcal{A}_M} \Bigg\vert \frac{1}{M} \sum_{m=1}^M \Big[ f(X_n^{(m)},A(X_n^{(m)}))  +\hat Y^{(m),A}_{n+1} \Big] - \mathbb{E}\Big[ J_n^{A,(\hat a_k^M)_{k=n+1}^{N-1}}(X_n) \Big] \Bigg\vert \Bigg]&  \nonumber\\
&\hspace{-10.5cm} = \int_{0}^{\infty} \mathbb{P}\Bigg[\sup_{A \in \mathcal{A}_M} \Bigg\vert \frac{1}{M} \sum_{m=1}^M \Big[ f(X_n^{(m)},A(X_n^{(m)}))  +\hat Y^{(m),A}_{n+1} \Big] - \mathbb{E}\Big[ J_n^{A,(\hat a_k^M)_{k=n+1}^{N-1}}(X_n) \Big] \Bigg\vert > \eps \Bigg] \diff \eps \nonumber\\
& \hspace{-10.5cm} = \int_{0}^{M'}  \mathbb{P}\Bigg[\sup_{A \in \mathcal{A}_M} \Bigg\vert \frac{1}{M} \sum_{m=1}^M \Big[ f(X_n^{(m)},A(X_n^{(m)}))  +\hat Y^{(m),A}_{n+1} \Big] - \mathbb{E}\Big[ J_n^{A,(\hat a_k^M)_{k=n+1}^{N-1}}(X_n) \Big] \Bigg\vert > \eps \Bigg] \diff \eps \nonumber \\
&\hspace{-10.3cm}  + \int_{M'}^{\infty}    \mathbb{P}\Bigg[\sup_{A \in \mathcal{A}_M} \Bigg\vert \frac{1}{M} \sum_{m=1}^M \Big[ f(X_n^{(m)},A(X_n^{(m)}))  +\hat Y^{(m),A}_{n+1} \Big] - \mathbb{E}\Big[ J_n^{A,(\hat a_k^M)_{k=n+1}^{N-1}}(X_n) \Big] \Bigg\vert > \eps \Bigg] \diff \eps \nonumber \\
& \hspace{-10.5cm} \leq \sqrt{2} \frac{(N-n) \lVert f \rVert_\infty + \lVert g \rVert_\infty}{\sqrt{M}} \nonumber \\
& \hspace{-10cm}+ 4\int_{0}^{\infty} \mathbb{P}\Bigg[\sup_{A \in \mathcal{A}_M} \Bigg\vert \frac{1}{M} \sum_{m=1}^M \Big[ f(X_n^{(m)},A(X_n^{(m)}))  +\hat Y^{(m),A}_{n+1} \nonumber\\ 
& \hspace{-6cm}- f(X'^{(m)}_n,A(X'^{(m)}_n))  -\hat Y^{\prime \, (m),A}_{n+1} \Big]   \Bigg\vert  > \eps\Bigg] \diff \eps.
\end{align}}
The second term in the r.h.s. of \eqref{ineq::interm333} comes from \eqref{ineq::step1}. It remains to write the latter as an expectation to obtain \reff{ineqstep2}.

\vspace{2mm}

\noindent \textit{Step 3: Introduction of additional randomness by random signs.}\\
Let $(r_m)_{1\leq m \leq M}$ be i.i.d. Rademacher r.v.\footnote{The probability mass function of a Rademacher r.v. is by definition $\frac{1}{2}\delta_{-1}+ \frac{1}{2}\delta_{1}$.}. We show that:
\begin{align}
\mathbb{E}\Bigg[\sup_{A \in \mathcal{A}_M} \Bigg\vert \frac{1}{M} \sum_{m=1}^M \Big[ f(X_n^{(m)},A(X_n^{(m)}))  +\hat Y^{(m),A}_{n+1} - f(X'^{(m)}_n,A(X'^{(m)}_n))  -\hat Y^{\prime \, (m),A}_{n+1} \Big]   \Bigg\vert \Bigg]  & \nonumber \\
&\hspace{-12cm} \leq 4\mathbb{E}\Bigg[\sup_{A \in \mathcal{A}_M} \Bigg\vert \frac{1}{M} \sum_{m=1}^M r_m \Big[ f(X_n^{(m)},A(X_n^{(m)}))  +\hat Y^{(m),A}_{n+1} \Big]\Bigg\vert \Bigg]. \label{ineq::speedLearning1}
\end{align}
Since for each $m=1,...,M$ the set of exogenous noises $(\eps'^{(m)}_k)_{n\leq k\leq N}$ and $(\eps^{(m)}_k)_{ n\leq k\leq N}$ are i.i.d., their  joint distribution remains the same if one 
randomly interchanges the correspon\-ding components. 
%$(\eps^{(1)}_k)_{n\leq k\leq N},...,(\eps^{(M)}_k)_{n\leq k\leq N},(\eps'^{(1)}_k)_{n\leq k\leq N},...,(\eps'^{(M)}_k)_{n\leq k\leq N}$ 
%remains the same if one randomly interchanges the corresponding components of  
%$\big((\eps^{(1)}_k)_{n\leq k\leq N},...,(\eps^{(M)}_k)_{n\leq k\leq N} \big)$ and $\big((\eps'^{(1)}_k)_{n\leq k\leq N},...,(\eps'^{(M)}_k)_{n\leq k\leq N}\big)$. 
Thus, the following holds for $\eps$ $\geq$ $0$:
{\small \begin{align}
\mathbb{P}\Bigg[\sup_{A \in \mathcal{A}_M} \Bigg\vert \frac{1}{M} \sum_{m=1}^M \Big[ f(X_n^{(m)},A(X_n^{(m)}))  +\hat Y^{(m),A}_{n+1} - f(X'^m_n,A(X'^m_n))  -\hat Y^{\prime \, (m),A}_{n+1} \Big]   \Bigg\vert > \eps \Bigg] & \nonumber \\
& \hspace{-12cm} = \mathbb{P}\Bigg[\sup_{A \in \mathcal{A}_M} \Bigg\vert \frac{1}{M} \sum_{m=1}^M r_m \Big[ f(X_n^{(m)},A(X_n^{(m)}))  +\hat Y^{(m),A}_{n+1} - f(X'^m_n,A(X'^m_n))  -\hat Y^{\prime \, (m),A}_{n+1} \Big]   \Bigg\vert > \eps \Bigg]  \nonumber \\
&\hspace{-12cm} \leq  \mathbb{P}\Bigg[\sup_{A \in \mathcal{A}_M} \Bigg\vert \frac{1}{M} \sum_{m=1}^M r_m  \Big[ f(X_n^{(m)},A(X_n^{(m)}))  +\hat Y^{(m),A}_{n+1} \Big]  \Bigg\vert > \frac{\eps}{2} \Bigg] \nonumber \\
& \hspace{-11.5cm} +  \mathbb{P}\Bigg[\sup_{A \in \mathcal{A}_M} \Bigg\vert \frac{1}{M} \sum_{m=1}^M r_m \Big[ f(X_n^{(m)},A(X_n^{(m)}))  +\hat Y^{(m),A}_{n+1} \Big]  \Bigg\vert > \frac{\eps}{2} \Bigg] \nonumber \\
&\hspace{-12cm} \leq  2  \mathbb{P}\Bigg[\sup_{A \in \mathcal{A}_M} \Bigg\vert \frac{1}{M} \sum_{m=1}^M r_m  \Big[ f(X_n^{(m)},A(X_n^{(m)}))  +\hat Y^{(m),A}_{n+1} \Big]  \Bigg\vert > \frac{\eps}{2} \Bigg]. \nonumber
\end{align}}
It remains to integrate on $\R_+$ w.r.t. $\eps$ to get \reff{ineq::speedLearning1}.

\vspace{2mm}

\noindent \textit{Step 4:} We show that 
\begin{align}
\mathbb{E}\Bigg[\sup_{A \in \mathcal{A}_M} \Bigg\vert \frac{1}{M} \sum_{m=1}^M r_m \Big[ f(X_n^{(m)},A(X_n^{(m)}))  +\hat Y^{(m),A}_{n+1} \Big]\Bigg\vert \Bigg] & \nonumber\\
& \hspace{-7cm} \leq  \frac{ (N-n) \lVert f \rVert_\infty +\lVert g \rVert_\infty }{\sqrt{M}} \nonumber\\
&\hspace{-6.5cm}+ \left( [f]_L +[f]_L  \sum_{k=n+1}^{N-1} \big(1+ \eta_M\gamma_M\big)^{k-n}[F]_L^{k-n} +  \eta_M^{N-n}\gamma_M^{N-n}[F]_L^{N-n}[g]_L \right)  \nonumber \\
& \hspace{-7cm}  = \mathcal{O}\bigg(  \frac{\gamma_M^{N-n}\eta_M^{N-n}}{\sqrt{M}}\bigg), \quad \text{ as $M\to \infty$}. \label{ineq::step3}
\end{align}
Adding and removing the cost obtained by control 0 at time $n$ yields:
{\small\begin{align}
\mathbb{E}\Bigg[\sup_{A \in \mathcal{A}_M} \Bigg\vert \frac{1}{M} \sum_{m=1}^M r_m \Big( f(X_n^{(m)},A(X_n^{(m)})) +\hat Y^{(m),A}_{n+1} \Big)\Bigg\vert \Bigg]  & \nonumber \\
&\hspace{-6.2cm} \leq  \mathbb{E}\Bigg[\bigg\vert \frac{1}{M} \sum_{m=1}^M r_m \Big( f(X_n^{(m)},0) + \hat Y^{(m),0}_{n+1}\Big)\bigg\vert \Bigg] \nonumber  \\
& \hspace{-6cm} +  \mathbb{E}\Bigg[\sup_{A \in \mathcal{A}_M} \Bigg\vert \frac{1}{M} \sum_{m=1}^M r_m \Big( f(X_n^{(m)},A(X_n^{(m)}))- f(X_n^{(m)},0)  + \hat Y^{(m),A}_{n+1} -  \hat Y^{(m),0}_{n+1}\Big)\Bigg\vert \Bigg]. \label{ineq::EEstep4}
\end{align}}
 We now bound the first term of the r.h.s. of \eqref{ineq::EEstep4}. 
By the Cauchy-Schwarz inequality, and recalling that $(r_m)_{1\leq m \leq M}$ are i.i.d. with zero mean such that $r_m^2=1$, we get
\beq
\mathbb{E}\Bigg[\bigg\vert \frac{1}{M} \sum_{m=1}^M r_m \Big( f(X_n^{(m)},0) + \hat Y^{(m),0}_{n+1} \Big)\bigg\vert \Bigg] &\leq & 
\frac{1}{M} \sqrt{\mathbb{E}\Bigg[\bigg\vert \sum_{m=1}^M r_m \Big( f(X_n^{(m)},0) +\hat Y^{(m),0}_{n+1} \Big)\bigg\vert^2 \Bigg]}  \nonumber \\
& \leq & \frac{1}{\sqrt{M}}  \big((N-n) \lVert f \rVert_\infty +\lVert g \rVert_\infty \big) \label{inter0step3}
\enq
Turn now to the second term of \eqref{ineq::EEstep4}.
By the Lipschitz continuity of $f$, it stands:
{\small \begin{align}
\mathbb{E}\Bigg[\sup_{A \in \mathcal{A}_M} \Bigg\vert \sum_{m=1}^M r_m \Big( f(X_n^{(m)},A(X_n^{(m)}))- f(X_n^{(m)},0)  + \hat Y^{(m),A}_{n+1}-  \hat Y^{(m),0}_{n+1} \Big)\Bigg\vert \Bigg]  & \nonumber \\
&\hspace{-9.4cm}  \leq [f]_L \E \bigg[ \sup_{A \in \mathcal{A}_M} \Big\vert  \sum_{m=1}^{M} r_{m}A(X_n^{(m)})   \Big\vert  \bigg] + \E \Bigg[ \sup_{A \in \mathcal{A}_M} \Bigg\vert \sum_{m=1}^{M} r_m \Big( \hat Y^{(m),A}_{n+1} - \hat Y^{(m),0}_{n+1}\Big)  \Bigg\vert  \Bigg]  \nonumber \\
& 	 \hspace{-9.4cm}   \leq \Bigg( [f]_L + [f]_L  \sum_{k=n+1}^{N-1} \big(1+ \eta_M\gamma_M\big)^{k-n} \E \left[ \Sup_{1 \leq m \leq M}\prod_{j=n+1}^{k}C\left(\eps_j^m\right)\right] \nonumber\\
& \hspace{-9cm} + [g]_L \left( 1+ \eta_M^{N-n}\gamma_M^{N-n} \right) \E \left[ \Sup_{1 \leq m \leq M}\prod_{j=n+1}^{N}C\left(\eps_j^m\right)\right] \Bigg) \E \left[ \sup_{A \in \mathcal{A}_M} \Big\vert  \sum_{m=1}^{M} r_{m}A(X_n^{(m)})   \Big\vert  \right] \label{interstep3bis}
\end{align}}
where we condition by the exogenous noise, use assumption \textbf{(HF-PI)} and the $\eta_M\gamma_M$-Lipschitz continuity of the estimated optimal controls at time k, for $k=n+1,\ldots,N-1$. Now, notice first that 
%$\E \Big[ \Sup_{1 \leq m \leq M}\prod_{j=n+1}^{N}C\big(\eps_j^m\big)\Big]$ can be bounded as follows:
\beq
\E \left[ \sup_{1 \leq m \leq M}  \prod_{k=n+1}^{N} C\left( \eps_k^m \right) \right] 
%&\leq \E \left[ \prod_{k=n+1}^{N}  \sup_{1 \leq m \leq M}  C\left( \eps_k^m \right) \right] 
&\leq &  \prod_{k=n+1}^{N} \E \left[  \sup_{1 \leq m \leq M}  C\left( \eps_k^m \right) \right] \; \leq \;  \rho_M^{N-n}, \label{inter1step3}
\enq
and moreover:
%It remains to bound $\E \bigg[ \Sup_{A \in \mathcal{A}_M} \Big\vert  \sum_{m=1}^{M} r_{m}A(X_n^{(m)})   \Big\vert  \bigg]$ as follows:
\beq
\E \bigg[ \sup_{A \in \mathcal{A}_M} \bigg\vert  \sum_{m=1}^{M} r_{m}A(X_n^{(m)})   \bigg\vert  \bigg] &\leq &  
\gamma_M \E \bigg[ \sup_{\lvert v \rvert_2 \leq 1/R}  \bigg\vert \sum_{m=1}^{M} r_m(v^T X_n^{(m)})_+  \bigg\vert \bigg] \nonumber \\
& \leq & \gamma_M \E \bigg[ \sup_{\lvert v \rvert_2 \leq 1/R}  \bigg\vert \sum_{m=1}^{M} r_mv^T X_n^{(m)} \bigg\vert \bigg],  \label{inter2step3} 
\enq
where $R>0$ is a bound for the state space (see e.g. the discussion on the Frank-Wolfe step p.10 of \cite{bach14} for a proof of this inequality), which implies by the Cauchy-Schwarz inequality: 
\beqs
 \E \bigg[ \sup_{A \in \mathcal{A}_M} \bigg\vert  \sum_{m=1}^{M} r_{m}A(X_n^{(m)})   \bigg\vert  \bigg] &\leq &  
 \frac{\gamma_M}{R} \sqrt{\E \bigg[   \bigg\vert \sum_{m=1}^{M} r_mX_n^{(m)}  \bigg\vert^2 \bigg]}  \leq    \gamma_M \sqrt{M} \label{ineq::EEstep4_2} 
 \enqs
since $(r_m)_m$ are i.i.d. Rademacher r.v..  Plug  \reff{inter1step3}   and \reff{inter2step3} into  \reff{interstep3bis} to obtain
{\small \begin{equation}
\label{ineq:LaterQ_truc}
\begin{split}
\mathbb{E}\Bigg[\sup_{A \in \mathcal{A}_M} \Bigg\vert \sum_{m=1}^M r_m \Big( f(X_n^{(m)},A(X_n^{(m)}))- f(X_n^{(m)},0)  + \hat Y^{(m),A}_{n+1}-  \hat Y^{(m),0}_{n+1} \Big)\Bigg\vert \Bigg] & \\
& 	 \hspace{-10cm}   \leq \Bigg( [f]_L + [f]_L  \sum_{k=n+1}^{N-1} \big(1+ \eta_M\gamma_M\big)^{k-n} \rho_M^{k-n} + [g]_L \left( 1+ \eta_M^{N-n}\gamma_M^{N-n} \right) \rho_M^{N-n} \Bigg) \gamma_M\sqrt{M}.
\end{split}
\end{equation}}
Plug then \eqref{inter0step3} and \eqref{ineq:LaterQ_truc}  into \eqref{ineq::EEstep4} to get \eqref{ineq::step3}.

\vspace{1mm}

\noindent \textit{Step 5: Conclusion.}  Plug \eqref{ineq::step3} into \eqref{ineq::speedLearning1} and combine it with \eqref{ineqstep2} to obtain the bound on the estimation error, as stated in 
\reff{estimesti} of Lemma \ref{speedLearning}. 
\ep

\section{Proof of Lemma \ref{lemapprox}}
\label{sec::annexerror1}

%We prove lemma \ref{lemapprox} in this section. Denote by $(a^{\mathrm{opt}}_n)_{0 \leq n\leq N-1}$ the optimal control solution of \eqref{defcontrol}. \\

Let $(\hat{a}^M_k)_{k=n+1}^{N-1}$  be the sequence of estimated controls at time $k=n+1,...,N-1$. 
Take $A \in \mathcal{A}_M$ and recall that we denote by $J^{A,( \hat{a})_{k=n+1}^{N-1}}_n$ the cost functional associated with the control $A$ at time $n$, and $\hat{a}^M_k$ at time $k=n+1, \ldots, N-1$. The latter is characterized as solution of the Bellman equation
\begin{equation*}
%\begin{cases}
\left\{
\begin{array}{ccl} 
J^{A,( \hat{a})_{k=n+1}^{N-1}}_N(x) &=& g(x) \\
J^{A,( \hat{a})_{k=n+1}^{N-1}}_n(x) &=& f(x, A(x)) + \mathbb{E}_{n,x}^A \big[J^{A,( \hat{a})_{k=n+1}^{N-1}}_{n+1}(X_{n+1}) \big],
\end{array}
%\end{cases}
\right.
\end{equation*}
where $\mathbb{E}_{n,x}^A[\cdot]$ stands for the expectation conditioned by $\{X_n=x \}$ when feedback control $A$ is followed at time $n$.

Take $n \in \{1,...,N\}$. The following holds:
\begin{align}
\eps^{\mathrm{approx}}_{\mathrm{PI},n} &:= \inf_{A \in \mathcal{A}_M}\E_M \big[ J^{A,( \hat{a})_{k=n+1}^{N-1}}_n(X_n) \big] 
- \inf_{A \in \A^{\Xc}}\E_M \big[ J^{A,( \hat{a})_{k=n+1}^{N-1}}_n(X_n) \big] \nonumber\\
& = \inf_{A \in \mathcal{A}_M}\E_M \big[ J^{A,( \hat{a})_{k=n+1}^{N-1}}_n(X_n) \big] - \E \big[  V_n(X_n) \big] \nonumber \\
& \hspace{.5cm}+\E\big[  V_n(X_n) \big] - \inf_{A \in \A^{\Xc}}\E_M \big[ J^{A,( \hat{a})_{k=n+1}^{N-1}}_n(X_n) \big] \nonumber \\
& \leq \inf_{A \in \mathcal{A}_M}\E_M \big[ J^{A,( \hat{a})_{k=n+1}^{N-1}}_n(X_n) \big] - \E \big[  V_n(X_n) \big], \label{ineq::approxError3}
\end{align}
where the last inequality stands because the value function is smaller than the cost functional associated with any other strategy.
We then apply the dynamic programming principle to obtain:
\begin{align}
	 \min_{A \in \mathcal{A}_M}\E_M \big[ J^{A,( \hat{a})_{k=n+1}^{N-1}}_n(X_n) \big] - \E \big[  V_n(X_n) \big] \nonumber \\
	  & \hspace{-4cm}\leq  \inf_{A \in \mathcal{A}_M}  \E_M \Big[ f\big( X_n, A(X_n)\big) + \E_n^A \big[ J^{( \hat{a})_{k=n+1}^{N-1}}_{n+1}(X_{n+1}) \big]  \Big] \nonumber \\
	& \hspace{-3.5cm} - \E \Big[  f\big(X_n,a^{\mathrm{opt}}_n(X_n)\big) +\E_n^{a^{\mathrm{opt}}} \big[V_{n+1}(X_{n+1}) \big] \Big]. \label{eq::boundApproxEstimation1} 
\end{align}
To bound the r.h.s. of \eqref{eq::boundApproxEstimation1}, first observe that for $A \in \mathcal{A}_M$:
\begin{align}
	\E_M \Big[ f\big( X_n, A(X_n)\big) + \E_n^A \big[ J^{( \hat{a})_{k=n+1}^{N-1}}_{n+1}(X_{n+1}) \big] \Big] & \nonumber\\ 
	&\hspace{-2cm}- \E \Big[  f\big(X_n,a^{\mathrm{opt}}_n(X_n)\big) +\E_n^{a^{\mathrm{opt}}}\big[ V_{n+1}(X_{n+1}) \big] \Big]   \nonumber\\
	& \hspace{-5.5cm}\leq \E \left[ \vert  f\big( X_n, A(X_n)\big) -f\big(X_n,a^{\mathrm{opt}}_n(X_n)\big) \vert \right] \nonumber \\
	& \hspace{-2cm}+ \E_M \Big[  \E_n^AJ^{( \hat{a})_{k=n+1}^{N-1}}_{n+1}(X_{n+1})-  \E_n^{a^{\mathrm{opt}}}V_{n+1}(X_{n+1}) \Big] \nonumber \\
	&\hspace{-5.5cm}\leq  \big( [f]_L  + \lVert V_{n+1} \rVert_\infty [r]_L \big) \E \left[ \vert  A(X_n) - a^{\mathrm{opt}}_n(X_n)\vert \right]  \nonumber\\
	& \hspace{-2cm}+ \lVert r \rVert_\infty  \E_M \Big[ J^{( \hat{a})_{k=n+1}^{N-1}}_{n+1}(X_{n+1})-  V_{n+1}(X_{n+1}) \Big],  \label{ineq::approxErrorEstimation2}
\end{align}
where we used twice assumption \ref{Hd} at the second-last line of \eqref{ineq::approxErrorEstimation2}. Inject 
\begin{align*}
	\E_M \Big[ J^{( \hat{a})_{k=n+1}^{N-1}}_{n+1}(X_{n+1}) \Big] &\leq \inf_{A \in \mathcal{A}_M} \E_M \Big[ J^{A,( \hat{a})_{k=n+2}^{N-1}}_{n+1}(X_{n+1}) \Big] + 2 \eps^{esti}_{n+1}
\end{align*}
into \eqref{ineq::approxErrorEstimation2} to obtain:
\begin{align}
\E_M \left[ f\big( X_n, A(X_n)\big) + \E_n^A \left[ J^{( \hat{a})_{k=n+1}^{N-1}}_{n+1}(X_{n+1}) \right] \right] &\nonumber\\
& \hspace{-2cm}- \E \Big[  f\big(X_n,a^{\mathrm{opt}}_n(X_n)\big) +\E_n^{a^{\mathrm{opt}}}\left[ V_{n+1}(X_{n+1}) \right] \Big]  &  \nonumber\\
& \hspace{-6cm}\leq   \big( [f]_L  + \lVert V_{n+1} \rVert_\infty [r]_L \big) \E \left[ \vert  A(X_n) - a^{\mathrm{opt}}_n(X_n)\vert \right] \nonumber \\
& \hspace{-5.5cm}+ \lVert r \rVert_\infty  \inf_{A \in \mathcal{A}_M} \E_M \left[ J^{A,( \hat{a})_{k=n+2}^{N-1}}_{n+1}(X_{n+1})- V_{n+1}(X_{n+1}) \right] + 2 \lVert r \rVert_\infty \eps^{esti}_{n+1}. \label{ineq::approxErrorEstimation3}
\end{align}
Plugging \eqref{ineq::approxErrorEstimation3} into \eqref{eq::boundApproxEstimation1} yields
\begin{align*}
	 \inf_{A \in \mathcal{A}_M}\E_M \left[ J^{A,( \hat{a})_{k=n+1}^{N-1}}_{n}(X_n) \right] - \E \left[  V_n(X_n) \right] & \\
	 & \hspace{-5cm}\leq  \lVert r \rVert_\infty  \inf_{A \in \mathcal{A}_M} \E_M \left[ J^{A,( \hat{a})_{k=n+2}^{N-1}}_{n+1}(X_{n+1})- V_{n+1}(X_{n+1}) \right] + 2 \lVert r \rVert_\infty \eps^{esti}_{n+1} \\
	 & \hspace{-4.5cm} +\big( [f]_L  +\lVert V_{n+1} \rVert_\infty [r]_L \big) \inf_{A \in \mathcal{A}_M}\E \left[ \vert  A(X_n) - a^{\mathrm{opt}}_n(X_n)\vert \right],
\end{align*}
which implies by induction, as $M\to +\infty$:
\begin{align*}
\E \left[ 	\inf_{A \in \mathcal{A}_M}\E_M \left[  J^{A,( \hat{a})_{k=n+1}^{N-1}}_{n}(X_n) \right] - \E \left[  V_n(X_n) \right] \right]  \\
& \hspace{-5.5cm} =\mathcal{O} \bigg(    \Sup_{n+1 \leq k \leq N-1} \E\left[ \eps^{esti}_{k} \right]+    \Sup_{n \leq k \leq N-1}\inf_{A \in \mathcal{A}_M}\E \left[ \vert  A(X_n) - a^{\mathrm{opt}}_n(X_n)\vert \right] \bigg).
\end{align*}
We now use Lemma \ref{speedLearning} to bound the expectations of the $\eps^{\mathrm{esti}}_{\mathrm{PI},k}$ for $ k=n+1, \ldots,N-1$, and plug the result into \eqref{ineq::approxError3} to complete the proof of Lemma \ref{lemapprox}.
\ep

\section{Function approximation by neural networks} \label{sec::functionApproximationUsingNeuralNetworks}

We assume $a^{\mathrm{opt}}_k \in \mathbb{L}^2(\mu)$, and show the relation  \eqref{eq:denseness1} in  Proposition \ref{ConvergenceError}.

\noindent The universal approximator theorem applies for $\mathcal{A}_\infty := \bigcup_{M=1}^\infty \mathcal{A}_M, $ and states that for all $\eps>0$, for all $k$, there exists a neural network $a^*$ in $\mathcal{A}_\infty$ such that:
$
\lVert a^{\mathrm{opt}}_k - a^* \rVert_{\infty} < \frac{\eps}{\mathcal{V}_d(\Xc)},
$
where $\mathcal{V}_d(\Xc)$ stands for the volume of compact set $\Xc$ seen as a compact of the Euclidean space $\mathbb{R}^d$. By integrating, we then get:
$
\sup_{n \leq k \leq N-1} \int_{\Xc}  \big\vert a^{\mathrm{opt}}_k(x) - a^*(x) \big \vert d \mu (x) < {\eps}.
$
Also, notice that $\big( \mathcal{A}_M \big)_{M \geq 1}$ is increasing, which implies that  $\mathcal{A}_\infty  = \lim_{M \to +\infty} \mathcal{A}_M$ and gives the existence of $M>0$ that depends on $\eps$ such that $a^* \in  \mathcal{A}_M$.

Therefore, we have shown that for $n=0,...,N-1$ 
\begin{equation*}
\sup_{n \leq k \leq N-1}\inf_{A \in \mathcal{A}_M} \E \big[ \vert A(X_k) - a^{\mathrm{opt}}_k(X_k) \vert  \big] \xrightarrow[M \to \infty]{} 0, \quad \text{ with }  X_k \sim \mu,
\end{equation*}
which is the required result stated in \eqref{eq:denseness1}.
\ep

\vspace{3mm}

\noindent We now show \eqref{eq:denseness2} of proposition \ref{ConvergenceError}. As stated in Proposition 6 of Section 4.7 in \cite{bach14}: we can approximate a $c$-Lipschitz function by a $\gamma_1$-norm less than $\gamma_M$ and uniform error less than $c \left( \frac{\gamma_M}{c}\right)^{-2d/(d+1)}\log \frac{\gamma_M}{c}$, and proposition 1 in \cite{bach14} shows that a function with $\gamma_1$ less than $\gamma_M$ may be approximated with $K_M$ neurons with uniform error $\gamma_M K_M^{-(d+3)/(2d)}$. \\
Thus, given $K_M$ and $ \gamma_M$, there exists a neural network $a^*$ in $\mathcal{V}_M$ such that
\begin{equation}
\label{ineq:bounderror_controle}
	\lVert a^* - a^{\mathrm{opt}} \rVert_{\infty} \leq c \left( \frac{\gamma_M}{c}\right)^{-2d/(d+1)}\log \left(\frac{\gamma_M}{c} \right) + \gamma_M K_M^{-(d+3)/(2d)}.
\end{equation}

\section{Proof of Lemma \ref{speedLearning_Hybrid}} \label{sec::estimationError_VI}

We prove Lemma \ref{speedLearning_Hybrid} in four steps. Since the proof is very similar to the one of Lemma \ref{speedLearning}, we only detail the arguments that are modified.

\noindent \textit{Step 1: Symmetrization by a ghost sample.} We take $\eps >0$ and show that for $M > 2\frac{\big( (N-n) \lVert f \rVert_\infty + \lVert g \rVert_\infty \big)^2}{\eps^2}$, the following holds
{\small
	\begin{align}
\mathbb{P}\Bigg[\sup_{A \in \mathcal{A}_M} \Bigg\vert \frac{1}{M} \sum_{m=1}^M \Big[ f(X_n^{(m)},A(X_n^{(m)}))  +\hat Y_{n+1}^{(m),A} \Big] - \mathbb{E}\Big[  f(X_n^{(m)},A(X_n^{(m)}))  +\hat Y_{n+1}^{(m),A} \Big] \Bigg\vert > \eps \Bigg]  \nonumber \\
& \hspace{-12.7cm}\leq 2 \mathbb{P}\Bigg[\sup_{A \in \mathcal{A}_M} \Bigg\vert \frac{1}{M} \sum_{m=1}^M\Big[ f(X_n^{(m)},A(X_n^{(m)}))  +\hat Y_{n+1}^{(m),A} -f(X'^{(m)}_n,A(X'^{(m)}_n))  -\hat Y^{\prime \, (m),A}_{n+1}\Big]  \Bigg\vert  > \frac{\eps}{2}\Bigg], \nonumber\\ \label{ineq::step1_hybrid}
\end{align}}
where $\big( X'^{(m)}_n \big)_{m=1}^M$ is an i.i.d. copy of $\big( X_n^{(m)} \big)_{m=1}^M$;  $\big( \eps'^m_{n+1} \big)_{m=1}^M$ is an i.i.d. copy of $\big( \eps^m_{n+1} \big)_{m=1}^M$; $
	\hat Y^{(m),A}_{n+1}:=  \hat{V}_{n+1}^M \left(F\Big( X_{n}^{(m)}, A \big(  X_{n}^{(m)} \big), \eps^m_{n+1} \Big)\right);$
	and finally
$
		\hat Y^{\prime \, (m),A}_{n+1}:=  \hat{V}_{n+1}^M \left(F\Big( X_{n}^{\prime (m)}, A \big(  X_{n}^{\prime (m)} \big), \eps'^m_{n+1} \Big)\right). \nonumber
$

Since $\hat{V}_{n}^M$ the estimated value function at time $n$, for $n$$=$$0,...,N-1$, is bounded by construction, \eqref{ineq::step1_hybrid} is proved the same as in step 1 of Lemma \ref{speedLearning}.
\ep

\vspace{3mm}
\noindent \textit{Step 2:}  The following result holds
{\small \begin{align}
\mathbb{E}\Bigg[\sup_{A \in \mathcal{A}_M} \Bigg\vert \frac{1}{M} \sum_{m=1}^M \Big[ f(X_n^{(m)},A(X_n^{(m)}))  +\hat Y^{(m),A}_{n+1} \Big] - \mathbb{E}\Big[  f(X_n^{(m)},A(X_n^{(m)}))  +\hat Y_{n+1}^{(m),A} \Big]  \Bigg\vert \Bigg] & \nonumber \\
&\hspace{-11.5cm} \leq 4\mathbb{E}\Bigg[\sup_{A \in \mathcal{A}_M} \Bigg\vert \frac{1}{M} \sum_{m=1}^M \Big[ f(X_n^{(m)},A(X_n^{(m)}))  +\hat Y^{(m),A}_{n+1} - f(X'^{(m)}_n,A(X'^{(m)}_n))  
-\hat Y^{\prime \, (m),A}_{n+1} \Big]   \Bigg\vert \Bigg] \nonumber \\
& \hspace{-11cm} + \mathcal{O}\left(  \frac{1}{\sqrt{M}}  \right), \label{step2HybridE}
\end{align}}
and it is proved the same way as step 2 in the proof of  Lemma \ref{speedLearning}.
\ep

\vspace{3mm}

\noindent \textit{Step 3: Introduction of additional randomness by random signs.} The following result holds:
\begin{align}
\mathbb{E}\Bigg[\sup_{A \in \mathcal{A}_M} \Bigg\vert \frac{1}{M} \sum_{m=1}^M \Big[ f(X_n^{(m)},A(X_n^{(m)}))  +\hat Y^{(m),A}_{n+1} - f(X'^{(m)}_n,A(X'^{(m)}_n))  -\hat Y^{\prime \, (m),A}_{n+1} \Big]   \Bigg\vert \Bigg]  & \nonumber \\
&\hspace{-12cm} \leq 4\mathbb{E}\Bigg[\sup_{A \in \mathcal{A}_M} \Bigg\vert \frac{1}{M} \sum_{m=1}^M r_m \Big[ f(X_n^{(m)},A(X_n^{(m)}))  +\hat Y^{(m),A}_{n+1} \Big]\Bigg\vert \Bigg], \label{ineq::speedLearning1_hybrid}
\end{align}
as proved already in step 3 in the proof of Lemma \ref{speedLearning}.
\ep

\vspace{3mm}

\noindent \textit{Step 4:} We show that 
\begin{align}
\mathbb{E}\Bigg[\sup_{A \in \mathcal{A}_M} \Bigg\vert \frac{1}{M} \sum_{m=1}^M r_m \Big[ f(X_n^{(m)},A(X_n^{(m)}))  +\hat Y^{(m),A}_{n+1} \Big]\Bigg\vert \Bigg] & \nonumber\\
& \hspace{-5cm}\leq  \frac{ (N-n) \lVert f \rVert_\infty +\lVert g \rVert_\infty}{\sqrt{M}} + \left( [f]_L + \rho_M\gamma_M \eta_M\right) \frac{\gamma_M}{\sqrt{M}}  \nonumber\\
& \hspace{-5cm}= \mathcal{O}\bigg(  \frac{\rho_M\gamma_M^2\eta_M}{\sqrt{M}}\bigg), \quad \text{ as $M\to +\infty$}. \label{ineq::step3_hybrid}
\end{align}

\noindent Adding and removing the cost obtained by control 0 at time $n$ yields:
\begin{align}
 \label{ineq::EEstep4_hybrid}
\mathbb{E}\Bigg[\sup_{A \in \mathcal{A}_M} \Bigg\vert \frac{1}{M} \sum_{m=1}^M r_m \Big( f(X_n^{(m)},A(X_n^{(m)})) +\hat Y^{(m),A}_{n+1} \Big)\Bigg\vert \Bigg]  & \\
& \hspace{-8cm}\leq  \mathbb{E}\Bigg[\bigg\vert \frac{1}{M} \sum_{m=1}^M r_m \Big( f(X_n^{(m)},0) + \hat Y^{(m),0}_{n+1}\Big)\bigg\vert \Bigg] \nonumber  \\
& \hspace{-7.6cm} +  \mathbb{E}\Bigg[\sup_{A \in \mathcal{A}_M} \Bigg\vert \frac{1}{M} \sum_{m=1}^M r_m \Big( f(X_n^{(m)},A(X_n^{(m)}))- f(X_n^{(m)},0)  + \hat Y^{(m),A}_{n+1} -  \hat Y^{(m),0}_{n+1}\Big)\Bigg\vert \Bigg]. \nonumber
\end{align}

\noindent The first term in the r.h.s. in \eqref{ineq::EEstep4_hybrid} is bounded as in the proof of Lemma \ref{speedLearning} by $\frac{(N-n)\lVert f \rVert_{\infty} +\lVert g \rVert_{\infty} }{\sqrt{M}}.$
We use the Lipschitz-continuity of $f$ as follows, to bound its second term:
\begin{align}
\mathbb{E}\Bigg[\sup_{A \in \mathcal{A}_M} \Bigg\vert \sum_{m=1}^M r_m \Big( f(X_n^{(m)},A(X_n^{(m)}))- f(X_n^{(m)},0)  + \hat Y^{(m),A}_{n+1}-  \hat Y^{(m),0}_{n+1} \Big)\Bigg\vert \Bigg]  & \nonumber \\
&\hspace{-10.5cm}  \leq [f]_L \E \bigg[ \sup_{A \in \mathcal{A}_M} \Big\vert  \sum_{m=1}^{M} r_{m}A(X_n^{(m)})   \Big\vert  \bigg] + \E \Bigg[ \sup_{A \in \mathcal{A}_M} \Bigg\vert \sum_{m=1}^{M} r_m \Big( \hat Y^{(m),A}_{n+1} - \hat Y^{(m),0}_{n+1}\Big)  \Bigg\vert \Bigg],  \nonumber \\ 
& 	 \hspace{-10.5cm}   \leq \left( [f]_L + \rho_M\eta_M\gamma_M\right) \E \Bigg[ \sup_{A \in \mathcal{A}_M} \Big\vert  \sum_{m=1}^{M} r_{m}A(X_n^{(m)})   \Big\vert  \Bigg], \nonumber
\end{align}
where we conditioned by the exogenous noise, used assumption \textbf{(HF)}, and  the $\eta_M\gamma_M$-Lipschitz continuity of the estimated value fonction at time $n+1$.

By using the same arguments as those presented to prove Lemma \ref{speedLearning}, we have $
\E \bigg[ \sup_{A \in \mathcal{A}_M} \bigg\vert  \sum_{m=1}^{M} r_{m}A(X_n^{(m)})   \bigg\vert  \bigg]	\leq \gamma_M \sqrt{M}, \label{ineq::EEstep4_2_hybrid} 
$
and then conclude that \eqref{ineq::step3_hybrid} holds.

\noindent \textit{Step 5: Conclusion}\\
Combining\eqref{step2HybridE},\eqref{ineq::speedLearning1_hybrid} and \eqref{ineq::step3_hybrid} results in the bound on the estimation error as stated in \eqref{estimesti_hybrid}.
\ep
\vspace{3mm}

\section{Proof of Lemma \ref{bound:approxErrorHybrid}} \label{sec:proofofLemmaApproxErrorHybrid}

We divide the proof of  Lemma \ref{bound:approxErrorHybrid} into two steps. First write
\begin{align}
\label{ineq:boundApproxErrorHybrid1}
\eps^{\mathrm{approx}}_{\mathrm{HN},n} &\leq  \inf_{A \in \mathcal{A}_M} \E_M\left[ f\left( X_n, A(X_n)\right) + \hat{V}_{n+1}^M \left( X_{n+1}^A \right) \right] - \E \left[ V_n(X_n)  \right] \nonumber\\
& \hspace{1cm} +  \E \left[ V_n(X_n)  \right]   - \inf_{A \in \A^{\Xc}} \E_M \left[ f\left( X_n, A(X_n)\right) + \hat{V}_{n+1}^M \left( X_{n+1}^A \right) \right]. 
\end{align}
\noindent \textit{Step 1:} We show
\begin{equation}
\begin{split}
\inf_{A \in \mathcal{A}_M} \E_M\left[ f\left( X_n, A(X_n)\right) + \hat{V}_{n+1}^M \left( X_{n+1}^A \right) \right] - \E \left[ V_n(X_n)  \right] & \\
&\hspace{-7cm}\leq \left( [f]_L + \lVert V_{n+1} \rVert_{\infty} [r]_L  \right)  \inf_{A \in \A^{\Xc}} \E_M \left[ \left\vert A(X_n) - a^{\mathrm{opt}}_n(X_n) \right\vert \right]  \\
& \hspace{-5.5cm}+  \lVert r \rVert_{\infty} \E_M \left[ \left\vert V_{n+1}(X_{n+1}) - \hat{V}_{n+1}^M(X_{n+1}) \right\vert \right]. 
\end{split}
 \label{ineq:Lem44_step4}
\end{equation}
Take $A \in \mathcal{A}_M$, and apply the dynamic programming principle to write 
{\small \begin{align}
\E_M\left[ f\left( X_n, A(X_n)\right) + \hat{V}_{n+1}^M \left( X_{n+1}^A \right) \right] - \E \left[ V_n(X_n)  \right] & \nonumber\\
& \hspace{-7cm}\leq  [f]_L \E_M \left[ \left\vert A(X_n) - a^{\mathrm{opt}}_n(X_n) \right\vert \right]  + \E_M \left[  \E_M \left[ \hat{V}_{n+1}^M(X_{n+1}^A)  \right] - \E_M \left[ V_{n+1} \left( X_{n+1}^{a^{\mathrm{opt}}_n}\right) \right]\right] \nonumber \\
& \hspace{-7cm}\leq  \left( [f]_L   +  \lVert V_{n+1} \rVert_{\infty} [r]_L  \right)\E_M\!\left[ \left\vert A(X_n) - a^{\mathrm{opt}}_n(X_n) \right\vert \right] +  \E_M\!\left[ \left\vert  \hat V_{n+1}^M (X_{n+1}^A) - V_{n+1} (X_{n+1}^A) \right\vert \right]\nonumber,
\end{align}}
where we used \ref{Hd} at the second-to-last line. By using one more time assumption \ref{Hd}, we then get:
\begin{align}
\E_M\left[ f\left( X_n, A(X_n)\right) + \hat{V}_{n+1}^M \left( X_{n+1}^A \right) \right] - \E \left[ V_n(X_n)  \right] & \nonumber\\
& \hspace{-5cm}\leq  \left( [f]_L   +  \lVert V_{n+1} \rVert_{\infty} [r]_L  \right)\E_M \left[ \left\vert A(X_n) - a^{\mathrm{opt}}_n(X_n) \right\vert \right] \nonumber \\
& \hspace{-4.5cm}+ \lVert r \rVert_{\infty} \E_M \left[ \left\vert  \hat V_{n+1}^M (X_{n+1}) - V_{n+1} (X_{n+1}) \right\vert \right], \quad \text{ with } X_{n+1} \sim \mu,\nonumber
\end{align}
which is the result stated in  \eqref{ineq:Lem44_step4}.

\vspace{3mm}
\noindent \textit{Step 2:} We show
\begin{equation}
\label{ineq:Lem44_step2}
\begin{split}
\E \left[ V_n(X_n)  \right]   - \inf_{A \in \A^{\Xc}} \E_M \left[ f\left( X_n, A(X_n)\right) + \hat{V}_{n+1}^M \left( X_{n+1}^A \right) \right]& \\
& \hspace{-5cm}\leq  \lVert r \rVert_{\infty} \E_M \left[ \left\vert V_{n+1}(X_{n+1}) - \hat{V}_{n+1}^M(X_{n+1}) \right\vert \right].
\end{split}
\end{equation}
Write
{\small \begin{align*}
 & \E\left[V_{n}(X_{n})\right]-\inf_{A\in\A^{\Xc}}\E_{M}\left[f\left(X_{n},A(X_{n})\right)+\hat{V}_{n+1}^{M}\left(X_{n+1}^{A}\right)\right]\\
 & \leq\inf_{A\in\A^{\Xc}}\E_{M}\left[f\left(X_{n},A(X_{n})\right)+{V}_{n+1}\left(X_{n+1}^{A}\right)\right]-\inf_{A\in\A^{\Xc}}\E_{M}\left[f\left(X_{n},A(X_{n})\right)+\hat{V}_{n+1}^{M}\left(X_{n+1}^{A}\right)\right]\\
 & \leq\inf_{A\in\A^{\Xc}}\E_{M}\left[V_{n+1}\left(X_{n+1}^{A}\right)-\hat{V}_{n+1}^{M}\left(X_{n+1}^{A}\right)\right]\leq\lVert r\rVert_{\infty}\E_{M}\left[\left\vert V_{n+1}(X_{n+1})-\hat{V}_{n+1}^{M}(X_{n+1})\right\vert \right]
\end{align*}}
which completes the proof of \eqref{ineq:Lem44_step2}. 

\vspace{3mm}

\noindent \textit{Step 3} \textit{Conclusion:}\\
We complete the proof of Lemma \ref{bound:approxErrorHybrid} by plugging \eqref{ineq:Lem44_step4} and \eqref{ineq:Lem44_step2} into  \eqref{ineq:boundApproxErrorHybrid1}. 
\ep

\vspace{3mm}

\section{Some useful Lemmas for the proof of Theorem \ref{theoKohler}}

\label{sectionKoh}

Fix $M$ $\in$ $\N^*$,  let $x_1,\ldots,x_M \in \R^d$, and  set  $x^M=(x_1,\ldots,x_M)$. Define the distance $d_2(f,g)$ between $f: \R^d \to \R $ and $g: \R^d \to \R$ by
\[
d_2(f,g)= \left(\frac{1}{M} \sum_{m=1}^{M} \left| f(x_m)-g(x_m) \right|^2 \right)^{1/2}. 
\]
An $\eps$-cover of $\mathcal{V}$ is a set of functions $f_1, \ldots, f_P : \R^d \to \R$ such that
\[
\min_{p =1,\ldots,P} d_2\left(f,f_p\right) < \eps, \quad \text{ for } f \in \mathcal{V}.
\]
Let $\mathcal{N}_2(\eps,\mathcal{V}, x^M)$  denote the size of the smallest $\eps$-cover of $\mathcal{V}$ w.r.t. the distance $d_2$, and set by convention $\mathcal{N}_2(\eps,\mathcal{V}, x^M)= + \infty$ if there does not exist any $\eps$-cover of $\mathcal{V}$ of finite size. $\mathcal{N}_2(\eps,\mathcal{V}, x^M)$ is called $\L^2$-$\eps$-covering number of $\mathcal{V}$ on $x^M$.
  
\vspace{3mm}

\begin{lemma} 
	\label{lem:H1} Let $(X,Y)$ be a random variable.
	Assume $\vert Y \vert \leq L $ a.s. and let $m(x)= \E[Y \vert X=x].$ 
	Assume $Y-m(X)$ is sub-Gaussian in the sense that\\
	$
	\Max_{m=1,...,M} c^2 \E \Big[ e^{(Y-m(X))^2/c^2}-1 \vert X \Big] \leq \sigma^2 \quad \text{ a.s.}
	$ for some $c,\sigma >0$. Let $\gamma_M, L \geq 1$ and assume that the regression function is bounded by $L$ and that $\gamma_M \xrightarrow[M \to +\infty]{}  +\infty$. Set \vspace{-2.5mm}
	\[
	\hat{m}_M:= \argmin_{\Phi \in \mathcal{V}_M} \frac{1}{M} \sum_{m=1}^{M} \big\vert \Phi(x_i)- Y_m \big\vert^2
	\] for some $\mathcal{V}_M$ of functions $\Phi : \R^d \to [-\gamma_M,\gamma_M]$ and some random variables ${Y}_1,...,{Y}_M$ which are bounded by $L$, and denote $$\Omega_g:= \bigg\{ f-g: f \in \mathcal{V}_M, \frac{1}{M} \sum_{m=1}^{M} \big\vert f(x_m)- g(x_m) \big\vert^2 \leq \frac{\delta }{\gamma_M^2} \bigg\}.$$ Then there exist constants $c_1,c_2>0$ which depend only on $\sigma$ and $c$ such that for any $\delta_M>0$ with \[
	\delta_M \xrightarrow[M \to +\infty]{} 0, \quad \frac{M \delta_M}{\gamma_M} \xrightarrow[M \to +\infty]{} +\infty, 
	\] 
	\begin{equation}
	\small
	\label{ineq:condtionCoveringNumber}
	c_1\frac{\sqrt{M} \delta }{\gamma_M^2}\geq  \int_{c_2\delta/\gamma_M^2}^{\sqrt{\delta}}  \log\left(  \mathcal{N}_2\bigg( \frac{u}{4\gamma_M}, \Omega_g, x^M_1 \bigg) \right)^{1/2} \diff u
	\end{equation} for all $\delta \geq \delta_M$ and all $g \in \mathcal{V}_M \cup \{m\}$ we have as $M \to +\infty$:
	\begin{align*}
	\E \Big[ \big\vert \hat{m}_M(X) - m(X) \big\vert^2 \Big] = \Op \left( \delta_M + \Inf_{ \Phi \in \mathcal{V}_M} \E \Big[ \big\vert \Phi(X) - m(X) \big\vert^2  \Big]  \right).
	\end{align*}
\end{lemma}
\begin{proof}
We refer to \cite{koh06} (see its Theorem 3) for a proof.
\end{proof}

\vspace{3mm}

\begin{lemma}
	\label{lem:H2}
	For any $\eps$ $>$ $0$, we have 
	\begin{equation*}
	\label{ineq:majoreCoveringNumber}
	\mathcal{N}_2 \Big( \eps, \mathcal{V}_M, (X_n^{(m)})_{1 \leq m \leq M} \Big) \leq \left( \frac{12 e \gamma_M \big( K_M +1\big)}{\eps} \right)^{(4d+9)K_M+1},
	\end{equation*}
	where the class of neural networks $\mathcal{V}_M$ is defined in Section \ref{subsectionHybrid}.
\end{lemma}
\begin{proof}
	We refer to \cite{kohkrztod10} for a proof.
\end{proof}

\newpage

\bibliography{bibtex}
\bibliographystyle{siamplain}

\end{document}